\documentclass[12pt]{article}
\usepackage{amssymb}

\usepackage{graphicx}
\usepackage{amsmath}

\input{tcilatex}
\input{tcilatex}

\begin{document}

\title{A Globally Accelerated Numerical Method For Optical Tomography With
Continuous Wave Source}
\author{Hua Shan$^{\bigtriangledown }$, Michael V. Klibanov$^{\ast }$, Jianzhong Su$%
^{\bigtriangledown },$\ Natee Pantong$^{\bigtriangledown },$ and Hanli Liu$%
^{\rhd }$ \\
$^{\bigtriangledown }$Department of Mathematics\\
University of Texas at Arlington\\
Arlington, TX 76019\\
$^{\ast }$Department of Mathematics and Statistics\\
University of North Carolina at Charlotte,\\
Charlotte, NC 28223\\
$^{\rhd }$ Department of Bioengineering\\
University of Texas at Arlington\\
Arlington, TX 76019\\
E-mail: mklibanv@uncc.edu\\
}
\maketitle

\begin{abstract}
A new numerical method for an inverse problem for an elliptic equation with
unknown potential is proposed. In this problem the point source is running
along a straight line and the source-dependent Dirichlet boundary condition
is measured as the data for the inverse problem. A rigorous convergence
analysis shows that this method converges globally, provided that the
so-called tail function is approximated well. This approximation is verified
in numerical experiments, so as the global convergence. Applications to
medical imaging, imaging of targets on battlefields and to electrical
impedance tomography are discussed.
\end{abstract}

\section{Introduction}

The phenomenon of multiple local minima and ravines of least squares
residual functionals represents the major obstacle for reliable numerical
solutions of Multidimensional Coefficient Inverse Problems (MCIPs) for
Partial Differential Equations (PDEs). We believe that because of the
applied nature of the discipline of Inverse Problems, the issue of
addressing the problem of local minima has \emph{vital importance} for this
discipline. Indeed, any gradient-like optimization method for such a
functional would likely converge to a local minimum, which is located far
from the correct solution. Furthermore, a global minimum, even a well
pronounced one, is not necesseraly located close to the true solution,
because of the ill-posed nature of MCIPs. Because of this, the vast majority
of current numerical methods for MCIPs are locally convergent ones, like,
for example Newton-like method, see, e.g., [1],[2],[4],[5] and many issues
of Inverse Problems. That is, convergence of such a method to the true
solution is rigorously guaranteed only if the initial guess is located
sufficiently close to that solution. However, in the majority of
applications such as e.g., medical and military ones, the media of interest
is highly heterogeneous, which means that a good first guess is unknown. The
latter naturally raises the question about the reliability of locally
convergent numerical methods for those applications, and this question is
well known to many practitioners working on computations of real world MCIPs.

Thus, we are interested in the issue of globally convergent numerical
methods for MCIPs. We call a numerical method \emph{globally convergent} if
the following two conditions are in place: (1) a rigorous convergence
analysis ensures that this method leads to a good approximation of the true
solution \emph{regardless} on the availability of a first good guess, and
(2) numerical experiments confirm the said convergence property.

In this paper we present an ``almost'' globally convergent method for an
MCIP for the equation 
\begin{equation}
\Delta _{\mathbf{x}}u-a\left( \mathbf{x}\right) u=-\delta \left( \mathbf{x}-%
\mathbf{x}_{0}\right) ,\mathbf{x=}\left( x,z\right) \in \mathbb{R}^{2}, 
\tag{1.1}
\end{equation}
\begin{equation}
\lim_{\left| \mathbf{x}\right| \rightarrow \infty }u\left( \mathbf{x},%
\mathbf{x}_{0}\right) =0.  \tag{1.2}
\end{equation}
Here $\mathbf{x}_{0}$ is the source position, and this position is running
along a line to generate the data for the inverse problem. We use the word
``almost'', because we rigorously prove global convergence only assuming
that we know a good approximation for the so-called ``tail function'', i.e.,
we assume that we know a good approximation of the second term of the
asymptotic behavior of the function $\ln \left[ u\left( \mathbf{x},\mathbf{x}%
_{0}\right) \right] $ for $\left| \mathbf{x}_{0}\right| \rightarrow \infty .$
Since this approximation is unknown analytically, we have decided to use a
heuristic iterative ``accelarator'' for convergence of tails and to confirm
the desired convergence numerically. Assuming that our tail function is
close to the correct one, we prove a global convergence result, which does
not rely on a good first guess for the solution. This is why we call our
method ``globally accelerated''. From the numerical standpoint, another
advantage for using the accelerator is that it gives us an approximation for
the tail, which seems to be rather close to the actual tail and we observe
this numerically. The only drawback is that we cannot establish this
rigorously.

We assume throughout this paper that the function $a\in C^{\alpha }\left( 
\mathbb{R}^{2}\right) ,a\geq const.>0$ and $a(\mathbf{x})=k^{2}=const.>0$
for $\mathbf{x}\in \mathbb{R}^{2}\diagdown \Omega ,$ where $\alpha \in
\left( 0,1\right) .$ The classic theory implies that there exists unique
solution of the problem (1.1), (1.2) such that $u\in C^{2+\alpha }\left(
\left| \mathbf{x}-\mathbf{x}_{0}\right| \geq \epsilon \right) ,\forall
\epsilon >0.$

The first generation of globally convergent numerical methods has started
from the so-called convexification algorithm [6]. The convexification is
using the projection technique with respect to all variables, except of one,
and a stable layer stripping procedure with respect to the latter variable.
\ While the work [6] is concerned with time/frequency dependent data, our
publication [8] is applying the convexification to the case of the running
source, which is the same as one in this paper. 
In the mathematics literature some other numerical techniques providing the global
convergence (see, e.g., [15-20])) are available. Their numerical 
implementations can be found in [21-24].

\textbf{\ }Recently the development of the second generation of globally
convergent numerical methods\ was initiated in [3]. The idea of [3] was
originated from our earlier publication [9].\ The concept of [3] also
overlaps in part with the scheme of the current paper. However, unlike our
current case, the time dependent data resulted from a single measurement are
considered in [3]. Laplace transform of either hyperbolic or parabolic
equation of [3] leads to the equation $\Delta w-s^{2}c(\mathbf{x})w-a(%
\mathbf{x})w=-\delta \left( \mathbf{x}-\mathbf{x}_{0}\right) ,$ where $s>0$
is the parameter of the Laplace transform and the source position is fixed.
Compared with our case, the main advantage of this equation is that the
asymptotic behavior of tails at $s\rightarrow \infty $ is known.\
Specifically, in the case when the coefficient $c(\mathbf{x})$ is unknown, $%
\lim_{s\rightarrow \infty }\left( \ln w/s^{2}\right) =0,$ and similarly when 
$a(x)$ is unknown. Ultimately, the knowledge of these limits enables one to
prove a global convergence theorem in [3] without a heuristic assumption. 
\textbf{\ }

We are interested in the extension \ of the idea of [3] to the case of the
running source instead of the changing time or frequency.\ In other words,
we consider almost the same inverse problem as one in [8].\ However, instead
of the convexification of [8] we develop an analogue of the method [3]. A
numerical method, similar with one of this publication, was published in our
early work [9]. However, the treatment of tails in section 4 of [9] was
different from one of our case, and that is why the global convergence
property was not observed in [9]. We also refer to subsection 5.4 of [4] for
another treatment of tails for a Newton-like locally convergent method for
an MCIP with frequency dependent data. We now explain the underlying reason
of our difficulties with the tail function from the physics standpoint. In
the case of the time dependent data for a hyperbolic equation [3] the tail
function is close to the so-called ``first arrival wave''. It \ is well
known that the first arrival signal is very informative one. However, it is
unclear what the first arrival signal is in our case of the elliptic
equation (1.1) with the running source.

We now formulate our inverse problem.

\textbf{Inverse Problem. }\emph{Denote }$\mathbf{x}=\left( x,z\right) .$%
\emph{\ Let }$\Omega \subset \mathbb{R}^{2}$\emph{\ be a bounded domain \
and }$\Gamma =\partial \Omega .$\emph{\ Let }$B$\emph{\ be a constant.} 
\emph{Suppose that in (1.1) }$\mathbf{x}_{0}=\left( s,B\right) \notin 
\overline{\Omega }.$\emph{\ Determine the coefficient }$a(x)$\emph{\ for }$%
x\in \Omega ,$\emph{\ assuming that the following function }$\varphi \left( 
\mathbf{x,}s\right) $\emph{\ is given } 
\begin{equation}
u\left( \mathbf{x},s\right) =\varphi \left( \mathbf{x,}s\right) ,\forall 
\mathbf{x}\in \Gamma ,\forall s\in \left[ \underline{s},\overline{s}\right] ,
\tag{1.3}
\end{equation}
\emph{where }$\overline{s}$\emph{\ is a sufficiently large number, }$%
\underline{s}<\overline{s}$\emph{\ is a certain fixed number and }$\left\{ 
\mathbf{x}=\left( s,B\right) ,s\geq \underline{s}\right\} \cap \partial
\Omega $\emph{=}$\varnothing .$

We consider the 2-D case for the sake of simplicity only for this
complicated problem.\ Generalizations of our method on the 3-D case are
feasible. The parameter count shows that the data $\varphi \left( \mathbf{x,}%
s\right) $ depends on two free parameters, so as the unknown coefficient $a(%
\mathbf{x}).$ Hence, this Inverse Problem is non-overdetermined. This
inverse problem has applications in imaging using light propagation in a
diffuse medium. This is the so-called continuous-wave (CW) light. In this
case the coefficient $a(\mathbf{x})$ is 
\begin{equation}
a(\mathbf{x})=3\left( \mu _{s}^{\prime }\mu _{a}\right) (\mathbf{x}), 
\tag{1.4}
\end{equation}
where $\mu _{s}^{\prime }(\mathbf{x})$ is the reduced scattering coefficient
and $\mu _{a}(\mathbf{x})$ is the absorption coefficient of the medium [1].
The first example of this application is imaging of targets on battlefields
covered by smog and flames using propagation of light originated by lasers.
In this application the laser source should be moved along a line and the
measurements of the output light should be performed at the boundary of the
domain of interest. Interestingly, the diffuse-like propagation of light
would be even helpful, because the direct light can miss the target. The
second applied example is in imaging of human organs or small animals using
near infrared light propagation. Note that this application is discussed in
many publications, in which locally convergent numerical methods are
developed, see, e.g., [1],[9],[13]. Also, the above Inverse Problem has
applications in Electrical Impedance Tomography, in which case the original
equation is $\nabla \cdot \left( \sigma \left( \mathbf{x}\right) \nabla
v\right) =-\delta \left( \mathbf{x}-\mathbf{x}_{0}\right) $ and the standard
change of variables $u=v\sqrt{\sigma }$ reduces this equation to (1.1)
assuming that $\sigma \left( \mathbf{x}\right) =1$ in a neighborhood of the
source position $\mathbf{x}_{0}.$ Here $\sigma \left( \mathbf{x}\right) \geq
const.>0$ is the electric conductivity of the medium. We note here that the potential function
$a (x, y)$ needs to remain positive for some of our mathematical arguments to work. The application to 
Electrical Impedance Tomography is limited under this constraint. 

\section{Nonlinear Integral Differential Equation}

Since the function $u$ is positive by the maximum principle,  we can
consider the function $v=\ln u.$ Since the source $\mathbf{x}_{0}=\left(
s,B\right) \notin \Omega $, we obtain the following equation from (1.1) 
\begin{equation}
\Delta v+\left| \nabla v\right| ^{2}=a(\mathbf{x})\text{ \ in }\Omega , 
\tag{2.1}
\end{equation}
\begin{equation}
v\left( \mathbf{x},s\right) =\varphi _{1}\left( \mathbf{x},s\right) ,\text{ }%
\forall \left( \mathbf{x},s\right) \in \Gamma \times \left( A,\overline{s}%
\right) ,\text{ }  \tag{2.2}
\end{equation}
where $\varphi _{1}=\ln \varphi .$ To eliminate the unknown coefficient $a(%
\mathbf{x)}$ from equation (2.1), differentiate it with respect to $s$ and
let 
\begin{equation}
q\left( \mathbf{x},s\right) =\partial _{s}v\left( \mathbf{x},s\right) . 
\tag{2.3}
\end{equation}
Then 
\begin{equation}
v\left( \mathbf{x},s\right) =-\int\limits_{s}^{\overline{s}}q\left( \mathbf{x%
},\tau \right) d\tau +T(\mathbf{x}),\mathbf{x}\in \Omega ,s\in \left( 
\underline{s},\overline{s}\right]  \tag{2.4}
\end{equation}

In (2.4) $T(\mathbf{x})$ is the so-called ``tail function''. The exact
expression for this function is of course $T(\mathbf{x})=v\left( \mathbf{x},%
\overline{s}\right) .$ We know only the first term of the asymptotic
expansion of the function $v\left( \mathbf{x},\overline{s}\right) $ at $%
\overline{s}\rightarrow \infty $ (below). As it was pointed out in
Introduction, if we would know the second term also, as it is the case of
the time dependent data of [3], then we would be better off approximating
the tail function. However, the absence of the knowledge of this term
significantly complicates the matter compared with [3]. Thus, we develop
below a heuristic iterative procedure of an iterative approximation of the
function $T(\mathbf{x}),$ with the aim of funding such an approximation $%
T_{appr}\left( \mathbf{x}\right) $ that $\nabla T_{appr}\left( \mathbf{x}%
\right) \approx \nabla T\left( \mathbf{x}\right) .$

We obtain from (2.1)-(2.4) 
\begin{equation}
\Delta q-2\nabla q\int\limits_{s}^{\overline{s}}\nabla qd\tau +2\nabla
q\nabla T=0,  \tag{2.5}
\end{equation}
\begin{equation}
q\left( \mathbf{x},s\right) =\psi \left( \mathbf{x},s\right) ,\forall \left( 
\mathbf{x},s\right) \in \Gamma \times \left( A,\overline{s}\right) , 
\tag{2.6}
\end{equation}
where 
\begin{equation*}
\psi \left( \mathbf{x},s\right) =\partial _{s}\ln \varphi \left( \mathbf{x}%
,s\right) .
\end{equation*}
The problem (2.5), (2.6) is nonlinear. In addition both functions $q$ and $T$
are unknown here. Now the main question is \emph{How to approximate well
both functions }$q$ \emph{and} $T$\emph{\ using (2.5), (2.6)?} Indeed, if we
approximate them well (in a certain sense, specified below), then the target
coefficient $a(\mathbf{x})$ would be reconstructed easily via backwards
calculations, see section 3. An equation similar with (2.5) was derived in
the convexification method [6],[8]. However the major difference between our
method and the convexification is in the numerical solution of the problem
(2.5), (2.6). Indeed, it is solution of this problem which represents the
major difficulty here.

\section{Layer Stripping With Respect to the Source Position}

An analogue of the nonlinear equation of this section for a different CIP,
in which the original PDE was either hyperbolic or parabolic was previously
derived in [3].

\subsection{Nonlinear equation}

We approximate the function $q\left( x,s\right) $ as a piecewise constant
function with respect to the pseudo frequency $s.$ That is, we assume that
there exists a partition 
\begin{equation*}
\underline{s}=s_{N}<s_{N-1}<...<s_{1}<s_{0}=\overline{s},s_{i-1}-s_{i}=h
\end{equation*}
of the interval $\left[ \underline{s},\overline{s}\right] $ with
sufficiently small grid step size $h$ such that 
\begin{equation*}
q\left( \mathbf{x},s\right) =q_{n}\left( \mathbf{x}\right) \text{ for }s\in %
\left[ s_{n},s_{n-1}\right) .
\end{equation*}
Hence 
\begin{equation}
\int\limits_{s}^{\overline{s}}\nabla q\left( \mathbf{x},\tau \right) d\tau
=\left( s_{n-1}-s\right) \nabla q_{n}\left( \mathbf{x}\right)
+h\sum\limits_{j=1}^{n-1}\nabla q_{j}\left( \mathbf{x}\right) ,s\in \left(
s_{n},s_{n-1}\right] .  \tag{3.1}
\end{equation}
We approximate the boundary condition (2.6) as a piecewise constant
function, 
\begin{equation}
q_{n}\left( \mathbf{x}\right) =\psi _{n}\left( \mathbf{x}\right) ,\mathbf{x}%
\in \partial \Omega ,  \tag{3.2}
\end{equation}
where 
\begin{equation}
\psi _{n}\left( \mathbf{x}\right) =\frac{1}{h}\int\limits_{s_{n}}^{s_{n-1}}%
\psi \left( \mathbf{x},s\right) ds.  \tag{3.3}
\end{equation}

Hence, for $s\in \left[ s_{n},s_{n-1}\right) $ equation (2.5) can be
rewritten as 
\begin{equation}
\widetilde{L}_{n}\left( q_{n}\right) :=\Delta q_{n}-2\nabla q_{n}\cdot
\left( h\sum\limits_{j=1}^{n-1}\nabla q_{j}-\nabla T\right) -2h\left( \nabla
q_{n}\right) ^{2}=0.  \tag{3.4}
\end{equation}
In sections 4 and 5 we address the question on how to solve equations (3.4)
for functions $q_{n}$ with the boundary conditions (3.2). \ 

\subsection{Reconstruction of the target coefficient}

Suppose that functions $\left\{ q_{i}\right\} _{i=1}^{n}$ are approximated
via solving problems (3.2), (3.4) and that the tail function is also
approximated. Then we reconstruct the target coefficient $a\left( \mathbf{x}%
\right) $ by backwards calculations as follows. First, we reconstruct the
function $v\left( \mathbf{x},s_{n}\right) $ by (2.4) as 
\begin{equation}
v\left( \mathbf{x},s_{n}\right) =-h\sum\limits_{j=1}^{n}q_{j}+T\left( 
\mathbf{x}\right) .  \tag{3.5}
\end{equation}
In principle we can reconstruct the target coefficient $a$ from (2.1).
However, it is unstable to take second derivatives. Hence, we first
reconstruct the function $u\left( \mathbf{x},s_{n}\right) $ as 
\begin{equation}
u\left( \mathbf{x},s_{n}\right) =\exp \left[ v\left( \mathbf{x},s_{n}\right) %
\right] .  \tag{3.6}
\end{equation}
Next, we use equation (1.1) in the weak form as 
\begin{equation}
-\int\limits_{\Omega }\nabla u\nabla \eta _{k}d\mathbf{x}=\int\limits_{%
\Omega }au\eta _{k}d\mathbf{x,}  \tag{3.7}
\end{equation}
where the test function $\eta _{k}\left( \mathbf{x}\right) ,k=1,...,K$ is a
quadratic finite element of a computational mesh with $\eta _{k}\left( 
\mathbf{x}\right) \mid _{\partial \Omega }=0.$ The number $K$ is finite and
depends on the mesh we choose. Equalities (3.7) lead to a linear algebraic
system which we solve. Then we obtain the function $\left( au\right) \left( 
\mathbf{x}\right) ,$ 
\begin{equation*}
\left( au\right) \left( \mathbf{x}\right) \approx
\sum\limits_{k=1}^{K}\alpha _{k}\eta _{k}\left( \mathbf{x}\right) .
\end{equation*}
Hence, 
\begin{equation}
a\left( \mathbf{x}\right) \approx \frac{1}{u\left( \mathbf{x},s_{n}\right) }%
\sum\limits_{k=1}^{K}\alpha _{k}\eta _{k}\left( \mathbf{x}\right) . 
\tag{3.8}
\end{equation}

\section{The Tail Function}

\bigskip We consider in this section two procedures for obtaining sequential
approximations for the tail function. First we find a first guess for the
tail function using the asymptotic behavior of the solution of the problem
(1.1), (1.2) as $\left| \mathbf{x}_{0}\right| \rightarrow \infty ,$ as well
as boundary measurements. Second, we describe an iterative procedure with
respect to tails. We call the combination of these two procedures
``accelerators'', because they help us to accelerate convergence of our
method. We stress that we cannot prove convergence of the second procedure.\
However, we have observed it in our numerical experiments. In our numerical
experiments we have worked only with a rectangular domain.\ Hence, we assume
in this section that 
\begin{equation*}
\left( x,z\right) \in \Omega :=\left\{ \left( x,z\right)
:x_{1}<x<x_{2},z_{0}<z<B\right\} .
\end{equation*}
However, we do not yet know how to address the issue of tails in the case of
an arbitrary convex domain $\Omega .$

\subsection{The first guess for tails}

First, we construct an approximation called ''asymptotic tail''. This is our
first accelerator. We consider the fundamental solution of the problem
(1.1), (1.2) for the case $a(x,z)\equiv k^{2}$. This solution is 
\begin{equation*}
u_{0}={\frac{1}{2\pi }}K_{0}(k{|}(x-s,z-B)|),
\end{equation*}
where $K_{0}(z)$ a modified Bessel function. It is well known that the
asymptotic behavior of this function is 
\begin{equation}
K_{0}(z)=\sqrt{\frac{\pi }{2\left| z\right| }}e{^{-k\left| z\right| }}(1+O({%
\frac{1}{\left| z\right| }})),|z|\rightarrow \infty .  \tag{4.1}
\end{equation}
Represent now solution of the problem (1.1), (1.2) as the solution of the
following integral equation 
\begin{equation}
u(x,z,s)={\frac{1}{2\pi }}K_{0}(k{|}(x-s,z-z_{m})|)  \tag{4.2}
\end{equation}
\begin{equation*}
-\frac{1}{2\pi }\int\limits_{\Omega }K_{0}\left( k\sqrt{\left( x-\xi \right)
^{2}+\left( z-\eta \right) ^{2}}\right) \left[ a\left( \xi ,\eta \right)
-k^{2}\right] u\left( \xi ,\eta ,s\right) d\xi d\eta .
\end{equation*}
Let $S\left( x,z,s\right) =|(x,z)-(s,z_{m})|$. The geometric meaning of $S$
is illustrated in Figure 1. Introducing the function 
\begin{equation*}
U(x,z,s)=2\sqrt{2\pi S}e{^{kS}}\cdot u(x,z,s)
\end{equation*}
and taking into account (4.1) and (4.2), we obtain that 
\begin{equation*}
u(x,z,s)=\sqrt{\frac{\pi }{2S}}e{^{-kS}}\left[ 1+\tilde{g}(x,z)+O\left( 
\frac{1}{S}\right) \right] ,S\rightarrow \infty .
\end{equation*}
The function $\tilde{g}(x,z)$ is unknown and is independent of $S$ as $%
S\rightarrow \infty $. Hence, we obtain for the function $v=\ln u$%
\begin{equation}
v(x,z,s)=-kS+\frac{1}{2}\ln (\frac{\pi }{2S})+g(x,z)+O\left( \frac{1}{S}%
\right) ,S\rightarrow \infty ,  \tag{4.3}
\end{equation}
where the unknown function $g(x,z)$ is derived from $\tilde{g}(x,z).$

We approximate the function $g(x,z)$ by two different methods and the final
answer is the average of two. The number of light sources $N=3$ is taken in
all our numerical experiments when we approximate this function. We start at 
$z=z_{0}$ where the boundary values are known. We decompose the boundary
values of $v$ into

\begin{equation}
v(x,z_{0},s_{j})=-kS_{j}+\frac{1}{2}\ln (\frac{\pi }{2S_{j}})+g_{j}(x,z_{0})
\tag{4.4}
\end{equation}
for j=1,2,3. Then we average to obtain 
\begin{equation}
g(x,z_{0},s)=\frac{1}{3}\sum\limits_{m=0}^{3}g_{j}(x,z_{0}),  \tag{4.5}
\end{equation}
Note that in (4.5) one should actually put ''$\approx $'' sign instead of
''=''.

However, the above procedure (4.4)-(4.5) gives us the value of the tail
functions $v\left( x,z,\overline{s}\right) $ only at $z:=z_{0},$ i.e., $%
v\left( x,z_{0},\overline{s}\right) .$ Equation (4.3) provides an
approximation for all $(x,z)\in \Omega $ if we simply set $%
g(x,z,s)=g(x,z_{0},s)$. In our numerical experiments we found that this is
insufficient. Hence, we use the measurement data from a different angle,
which enhances our numerical results. We obtain a similar tail function
using the measurement data at the lower edge of $\Omega $, i.e., at $x=x_{1}$
and got a second tail function using the idea similar with the above.\ Thus,
we have approximated $v\left( x_{1},z,\overline{s}\right) .$ Finally we set
for the first guess for producing a tail function 
\begin{equation}
T_{1, 0}\left( x,z,\overline{s}\right) :=\frac{1}{2}\left[ v\left( x,z_{0},%
\overline{s}\right) +v\left( x_{1},z,\overline{s}\right) \right].  \tag{4.6}
\end{equation}

\subsection{The second accelerator: iterations with respect to tails}

The second accelerator involves another iterative process that enhances the
reconstructed inclusion. Recall that $k^{2}:=a_{0}$ is the constant
background outside of our domain $\Omega ^{\prime }$. We now show how to
find an approximation $T_{1}\left( x,z,\overline{s}\right) $ for the tail
function. Let $u_{1,0}=e^{v_{0}}$ where $v_{0}=T_{1,0}\left( x,z,\overline{s}%
\right) $ is the function introduced above. We reconstruct the approximation 
$a_{1,1}(x,z)$ for the unknown coefficient $a(x,z)$ using the tail function
(4.6) through the inversion formula in equation (3.7) for all quadratic
finite element $\eta _{k}$ : 
\begin{equation*}
-\int\limits_{\Omega }\nabla u_{1,0}\nabla \eta _{k}d\mathbf{x}%
=\int\limits_{\Omega }a_{1,1}u_{1,0}\eta _{k}d\mathbf{x.}
\end{equation*}
Next, we apply (3.8). Then on the second step we solve the following
boundary value problem 
\begin{equation*}
\Delta u_{1,1}-a_{1,1}\left( x,z\right) u_{1,1}=0,\left( x,z\right) \in
\Omega ,
\end{equation*}
\begin{equation*}
u_{1,1}\mid _{\partial \Omega }=\varphi \left( \mathbf{x},\overline{s}%
\right) .
\end{equation*}
The reason for doing so is that we need to satisfy the boundary condition
obtained from measurements. 

We now describe a heuristic idea which motivates our iterative scheme. Let
the function $u$ be the solution of the following boundary value problem 
\begin{equation*}
\Delta u-a\left( x,z\right) u=-\delta (x,z),\left( x,z\right) \in \Omega ,
\end{equation*}
\begin{equation*}
u\mid _{\partial \Omega }=\varphi \left( \mathbf{x},\overline{s}\right)
\end{equation*}
with the unknown coefficient $a\left( x,z\right) $ and the function $u_{0}$
satisfies 
\begin{equation*}
\Delta u_{0}-a_{0}u_{0}=-\delta (x,z),\left( x,z\right) \in \Omega ,
\end{equation*}
\begin{equation*}
u_{0}\mid _{\partial \Omega }=\varphi \left( \mathbf{x},\overline{s}\right)
\end{equation*}
with the background function $a_{0}=k^{2}.$ Denote $p=u-u_{0}.$ Then 
\begin{equation*}
\Delta p-a\left( x,z\right) p=\left( a\left( x,z\right) -a_{0}\right) u_{0},
\end{equation*}
\begin{equation*}
p\mid _{\partial \Omega }=0.
\end{equation*}

Motivated by this idea, we introduce an iterative scheme and repeat the
procedure until it converges. Suppose that after $m-1$ iterations we have
constructed the function $u_{1,m-1}$ and have found the approximation $%
a_{1,m}(x,z)>0,m\geq 1$ for the unknown coefficient $a(x,z)$ using equation
(3.7)-(3.8). Then on the iteration number $m$, we solve the following
boundary value problem: 
\begin{equation*}
\Delta p_{1,m}-a_{1,m}\left( x,z\right) p_{1,m}=\left( a_{1,m}\left(
x,z\right) -a_{1,m-1}\left( x,z\right) \right) u_{1,m-1},
\end{equation*}
\begin{equation*}
p_{1,m}\mid _{\partial \Omega }=0.
\end{equation*}
Next, we set 
\begin{equation*}
u_{1,m}=u_{1,m-1}+p_{1,m}.
\end{equation*}

To accelerate convergence, we modify the iterative scheme slightly to solve
the following boundary value problems: 
\begin{equation*}
\Delta p_{1,m}-a_{1,m}\left( x,z\right) p_{1,m}=\lambda _{m}\left(
a_{1,m}\left( x,z\right) -a_{1,m-1}\left( x,z\right) \right) u_{1,m-1}
\end{equation*}
where 
\begin{equation*}
\lambda _{m}={\frac{exp\{\pi ^{2}e^{-(m-1)}\left( a_{1,m}\left( x,z\right)
-a_{1,m-1}\left( x,z\right) \right) ^{2}\}}{{\gamma }^{m}}}
\end{equation*}
and $\gamma =1.05.$ This choice of $\lambda _{m}$ is made in numerical
experiments. The choice of $\lambda _{m}$ makes the sequence converge after
about 50 iterations, instead of more than 300 in cases where $\lambda
_{m}\equiv 1$.

Once we have $u_{1,m}=u_{1,m-1}+p_{1,m}$, we construct $a_{1,m+1}$ by
equation (3.7) in the form of 
\begin{equation*}
-\int\limits_{\Omega }\nabla u_{1,m}\nabla \eta _{l}d\mathbf{x}%
=\int\limits_{\Omega }a_{1,m+1}u_{1,m}\eta _{l}d\mathbf{x,}
\end{equation*}
for all quadratic finite elements $\eta _{l},l=1,,,,.,K$ and use (3.8) then.
We iterate until the process converges, i.e., 
\begin{equation*}
\frac{\left\| a_{1,m_{1}}-a_{1,m_{1}-1}\right\| _{L_{2}\left( \Omega \right)
}}{\left\| a_{1,m_{1}-1}\right\| _{C\left( \overline{\Omega }\right) }}\leq
\varepsilon 
\end{equation*}
for a small $\varepsilon >0$ of our choice, see (7.1) for a detail. We set
for the first approximation for the unknown coefficient 
\begin{equation*}
a_{1}\left( x,z\right) :=a_{1,m_{1}}(x,z).
\end{equation*}

Then we set for the tail 
\begin{equation}
T_{1}\left( x,z,\overline{s}\right) =\ln u_{1,m_{1}}\left( x,z\right) , 
\tag{4.7}
\end{equation}
assuming that $u_{1,m_{1}}>0.$ Then we proceed with calculating the
functions $q_{n}$ as in section 5.

Remarks 4.1. 1. Unfortunately we cannot yet prove that functions $a_{1,m}>0$%
. Therefore, we cannot prove analytically neither the existence of solutions
of the above Dirichlet boundary value problems for functions $p_{1,m}$ \ nor
the positivity of functions $u_{1,m}.$\ Neither we cannot analytically prove
that functions $u_{1,m}$
converge, nor that our tail $T_{1}$ is close to the correct tail $T.$
Nevertheless, we observe all these ''nice'' properties in our computations.
Figure 2 displays the comparison of graphs of tails side by side, they have
little visible difference. 2. Unlike [3], we do not change tails in all
subsequent steps when calculating functions $q_{n}.$ In other words, the
tail function is kept the same $T:=T_{1}\left( x,z,\overline{s}\right) $ in
all follow up steps of our algorithm.

\section{The Algorithm for Approximating Functions $q_{n}$}

\bigskip \emph{Step} $1$. We need to find an approximation for the function $%
q_{1}.$ To do this, we solve equation (3.4) with the boundary condition
(3.2) at $n=1$ iteratively for $q_{1}.$ That is, we should solve 
\begin{equation}
\Delta q_{1}+2\nabla q_{1}\nabla T_{1}=2h\left( \nabla q_{1}\right) ^{2} 
\tag{5.1}
\end{equation}
\begin{equation}
q_{1}\left( \mathbf{x}\right) =\overline{\psi }_{1}\left( \mathbf{x}\right) ,%
\mathbf{x}\in \partial \Omega ,  \tag{5.2}
\end{equation}
We solve the problem (5.1), (5.2) iteratively as 
\begin{equation}
\Delta q_{1,k}+2\nabla q_{1,k}\nabla T_{1}-2h\nabla q_{1,k}\nabla
q_{1,k-1}=0,q_{1,k}\left( \mathbf{x}\right) =\overline{\psi }_{1}\left( 
\mathbf{x}\right) ,\mathbf{x}\in \partial \Omega  \tag{5.3}
\end{equation}
where $q_{1,0}=0$.

We proceed with calculating the function $q_{1,m+1}$ as in (5.3). We iterate
in (5.3) until the process converges, i.e., 
\begin{equation*}
\left\| q_{1,k_{1}}-q_{1,k_{1}-1}\right\| _{L_{2}\left( \Omega \right) }\leq
\varepsilon
\end{equation*}
for a small $\varepsilon >0$ of our choice, and $\varepsilon $ is the same
as in (4.13)$.$ We set $q_{1}:=q_{1,k_{1}}$. The next reconstruction $%
a_{2}(x,z)$ is obtained using equations (3.5)-(3.8), where $T:=T_{1}$.

\emph{Step} $n$. We now find an approximation for the function $q_{n}$
assuming that functions $q_{1},...,q_{n-1}$ are found. We solve iteratively
equation (3.4) with the boundary condition (3.2) as follows 
\begin{equation}
\Delta q_{n,k}-2h\sum\limits_{j=1}^{n-1}\nabla q_{j}\cdot \nabla
q_{n,k}+2\nabla q_{n,k}\nabla T_{1}-2h\nabla q_{n,k}\nabla
q_{n,k-1}=0,k=1,...,m_{n},  \tag{5.4}
\end{equation}
\begin{equation}
q_{n,k}\left( \mathbf{x}\right) =\psi _{n}\left( \mathbf{x}\right) ,\mathbf{x%
}\in \partial \Omega ,  \tag{5.5}
\end{equation}
where $q_{n,0}:=q_{n-1}.$We iterate until the process converges, i.e., until 
\begin{equation*}
\left\| q_{n,k_{n}}-q_{n,k_{n}-1}\right\| _{L_{2}\left( \Omega \right) }\leq
\varepsilon
\end{equation*}
for the above small $\varepsilon >0.$ We set $q_{n}:=q_{n,k_{n}}$. Then $%
a_{n+1}(x,z)$ is reconstructed using equations (3.5)-(3.8), $T:=T_{1}$.

\bigskip We find functions $a_{1},...,a_{N}$ for $n=1,...,N,$ where $N$ is
the number of subintervals of the interval $\left[ s_{0},\overline{s}\right]
.$ Finally, the resulting function $a\left( x,z\right) $ is 
\begin{equation}
a\left( x,z\right) =\frac{1}{N}\sum\limits_{i=1}^{N}a_{i}\left( x,z\right) .
\tag{5.6}
\end{equation}

We stress that we did not prove convergence of tails $T_n$ nor $q_n$
rigorously.\ Neither we cannot prove that functions $p_{n,m}$ in section 4
are positive, because we cannot prove that $a_{n,m}-a_{n,m-1}<0$ (in order
to apply the maximum principle). However, we have observed both the
positivity of functions $p_{n,m}$ and convergence of tails $T_n$ and $q_n$
in our computations.

\section{Convergence}

Below we follow the concept of Tikhonov for ill-posed problems [10], which
is one of backbones of this theory. By this concept one should assume first
that there exists an ``ideal'' exact solution of the problem with the exact
data. Next, one should assume the presence of an error in the data of the
level $\zeta ,$ where $\zeta >0$ is a small parameter. Suppose that an
approximate solution is constructed for each sufficiently small $\zeta $.
This solution is called a ``regularized solution'', if the $\zeta -$%
dependent family of these solutions tends to that exact solution as $%
\varsigma $ tends to zero. Hence, one should prove this convergence (Theorem
6.1).

In this section we use the Schauder's theorem to estimate functions $q_{n.k},
$ see \S 1 of Chapter 3 of [7] for this theorem. Since the Schauder's
theorem requires $C^{2+\alpha }$ smoothness of the boundary $\partial \Omega 
$, we assume in this section that $\Omega \subset \mathbb{R}^{2}$ is a
convex bounded domain with $\partial \Omega \in C^{2+\alpha }.$ This is, of
course in a disagreement with the above case of $\Omega $ being a rectangle.
However, we use the rectangle only because of the problem with tails, in
which we cannot rigorously prove that they are small and do not yet know how
to approximate them well heruistically for the case of a more general
domain, so that they would be close to correct tails. However, an analogue
of our convergence result (Theorem 6.1) can be proven for the case when $%
\Omega $ is rectangle and an FEM (i.e., discrete) version of equation (3.4)
is considered with a fixed number $R$ of finite elements. To do this, one
needs to consider the weak formulation of (3.4) and to use the Lax-Milgram
theorem instead of the Schauder's theorem.\ Although the Lax-Milgram theorem
would provide only estimates of $H^{1}$ norms of functions $q_{n}$ rather
than more desirable $C^{2}$ norms, but using the equivalency of norms in
finite dimensional spaces, we can still get estimates of $C^{2}$ norms and
these estimates would naturally depend on $R$.

\subsection{Exact solution}

Following the Tikhonov concept, we need to introduce the definitions of the
exact solution first. We assume that there exists an exact coefficient
function $a^{\ast }\left( x\right) \in C^{\alpha }\left( \overline{\Omega }%
\right) ,$ $\alpha =const.\in \left( 0,1\right) ,$ which is a solution of
our Inverse Problem. Let the function 
\begin{equation*}
u^{\ast }\left( \mathbf{x},s\right) \in C^{2+\alpha }\left( \left| \mathbf{x}%
-\mathbf{x}_{0}\right| \geq \varepsilon \right) ,\forall \varepsilon
>0,\forall \mathbf{x}_{0}=\left( s,B\right) >0,\forall s\in \left[ 
\underline{s},\overline{s}\right]
\end{equation*}
be the solution of the problem (1.1), (1.2) with $a\left( \mathbf{x}\right)
:=a^{\ast }\left( \mathbf{x}\right) $. Let 
\begin{equation*}
v^{\ast }\left( \mathbf{x},s\right) =\ln u^{\ast }\left( \mathbf{x},s\right)
,q^{\ast }\left( \mathbf{x},s\right) =\frac{\partial v^{\ast }\left( \mathbf{%
x},s\right) }{\partial s},T^{\ast }\left( \mathbf{x},\overline{s}\right)
=v^{\ast }\left( \mathbf{x},\overline{s}\right) .
\end{equation*}
By (2.1) 
\begin{equation}
a^{\ast }\left( \mathbf{x}\right) =\Delta v^{\ast }+\left( \nabla v^{\ast
}\right) ^{2}.  \tag{6.1}
\end{equation}
Also, the function $q^{\ast }$ satisfies the following analogue of equation
(2.5) 
\begin{equation}
\Delta q^{\ast }-2\nabla q^{\ast }\cdot \int\limits_{s}^{\overline{s}}\nabla
q^{\ast }\left( x,\tau \right) d\tau +2\nabla q^{\ast }\nabla T^{\ast }=0 
\tag{6.2}
\end{equation}
with the boundary condition (see (2.6)) 
\begin{equation}
q^{\ast }\left( \mathbf{x},s\right) =\psi ^{\ast }\left( \mathbf{x},s\right)
,\left( \mathbf{x},s\right) \in \partial \Omega \times \left[ \underline{s},%
\overline{s}\right] ,  \tag{6.3}
\end{equation}
where $\psi ^{\ast }\left( \mathbf{x},s\right) =\partial _{s}\ln \varphi
^{\ast }\left( \mathbf{x},s\right) ,$ where $\varphi ^{\ast }\left( \mathbf{x%
},s\right) =u^{\ast }\left( \mathbf{x},s\right) $ for $\left( \mathbf{x}%
,s\right) \in \partial \Omega \times \left[ \underline{s},\overline{s}\right]
.$

\textbf{Definition}. We call the function $q^{\ast }\left( \mathbf{x}%
,s\right) $ \emph{the exact solution} of the problem (2.5), (2.6) with the 
\emph{exact boundary} condition $\psi ^{\ast }\left( x,s\right) $.
Naturally, the function $a^{\ast }\left( \mathbf{x}\right) $ from (6.1) is
called the \emph{exact solution} of our Inverse Problem.

Therefore, 
\begin{equation}
q^{\ast }\left( \mathbf{x},s\right) \in C^{2+\alpha }\left( \overline{\Omega 
}\right) \times C^{1}\left[ \underline{s},\overline{s}\right] .  \tag{6.4}
\end{equation}
We now approximate the function $q^{\ast }\left( x,s\right) $ via a
piecewise constant function with respect to $s\in \left[ \underline{s},%
\overline{s}\right] .$ Let 
\begin{equation*}
q_{n}^{\ast }\left( \mathbf{x}\right) =\frac{1}{h}\int%
\limits_{s_{n}}^{s_{n-1}}q^{\ast }\left( \mathbf{x},s\right) ds,\text{ }\psi
_{n}^{\ast }\left( \mathbf{x}\right) =\frac{1}{h}\int%
\limits_{s_{n}}^{s_{n-1}}\psi ^{\ast }\left( x,s\right) ds
\end{equation*}
Then by (6.4) 
\begin{equation}
q^{\ast }\left( x,s\right) =q_{n}^{\ast }\left( x\right) +Q_{n}\left(
x,s\right) ,\psi ^{\ast }\left( x,s\right) =\psi _{n}^{\ast }\left( x\right)
+\Psi _{n}\left( x,s\right) ,  \tag{6.5}
\end{equation}
$s\in \left[ s_{n},s_{n-1}\right] ,$ where functions $Q_{n},\Psi _{n}$ are
such that for $s\in \left[ s_{n},s_{n-1}\right] $ 
\begin{equation}
\left\| Q_{n}\left( \mathbf{x},s\right) \right\| _{C^{2+\alpha }\left( 
\overline{\Omega }\right) }\leq C^{\ast }h,\left\| \Psi _{n}\left( \mathbf{x}%
,s\right) \right\| _{C^{2+\alpha }\left( \overline{\Omega }\right) }\leq
C^{\ast }h,\forall s\in \left[ s_{n},s_{n-1}\right] ,n=1,...,N,  \tag{6.6}
\end{equation}
where the constant $C^{\ast }>0$ depends only on $C^{2+\alpha }\left( 
\overline{\Omega }\right) \times C^{1}\left[ \underline{s},\overline{s}%
\right] $ and $C^{2+\alpha }\left( \partial \Omega \right) \times C^{1}\left[
\underline{s},\overline{s}\right] $ norms of functions $q^{\ast }$ and $\psi
^{\ast }$ respectively. Hence 
\begin{equation}
q_{n}^{\ast }\left( \mathbf{x}\right) =\psi _{n}^{\ast }\left( \mathbf{x}%
\right) ,\mathbf{x}\in \partial \Omega  \tag{6.7}
\end{equation}
and the following analogue of equation (3.4) holds 
\begin{equation}
\widetilde{L}_{n}\left( q_{n}\right) :=\Delta q_{n}^{\ast }-2\nabla
q_{n}^{\ast }\cdot \left( h\sum\limits_{j=1}^{n-1}\nabla q_{j}^{\ast
}-\nabla T^{\ast }\right) -2h\left( \nabla q_{n}^{\ast }\right)
^{2}=F_{n}^{\ast }\left( \mathbf{x},h\right) .  \tag{6.8}
\end{equation}
where the function $F_{n}\left( \mathbf{x},h\right) \in C^{\alpha }\left( 
\overline{\Omega }\right) $ and 
\begin{equation}
\left\| F_{n}^{\ast }\left( \mathbf{x},h\right) \right\| _{C^{\alpha }\left( 
\overline{\Omega }\right) }\leq C^{\ast }h.  \tag{6.9}
\end{equation}

We also assume that the data $\varphi \left( \mathbf{x},s\right) $ in (1.3)
are given with an error. This naturally produces an error in the function $%
\psi \left( \mathbf{x},s\right) $ in (2.6). An additional error is
introduced due to taking the average value of $\psi \left( \mathbf{x}%
,s\right) $ over the interval $\left( s_{n},s_{n-1}\right) $. Hence, it is
reasonable to assume that 
\begin{equation}
\left\| \overline{\psi }_{n}^{\ast }\left( \mathbf{x}\right) -\psi
_{n}\left( \mathbf{x}\right) \right\| _{C^{2+\alpha }\left( \partial \Omega
\right) }\leq C_{1}\left( \sigma +h\right) ,  \tag{6.10}
\end{equation}
where $\sigma >0$ is a small parameter characterizing the level of the error
in the data $\psi \left( \mathbf{x},s\right) $ and the constant $C_{1}>0$ is
independent on numbers $\sigma $, $h$ and $n$.

\textbf{Remark 6.1.} It should be noted that usually the data $\varphi
\left( \mathbf{x},s\right) $ in (1.3) are given with a random noise.
Although the differentiation of the noisy data is an ill-posed problem, but
there exist effective numerical regularization methods of its solution. We
are not addressing the corresponding theory here referring the reader to
e.g., [4], and also see section 7 for our way of handling it.

\subsection{Convergence theorem}

First, we reformulate the Schauder's theorem in a way, which is convenient
for our case. Introduce the positive constant $M^{\ast }$ as 
\begin{equation*}
M^{\ast }=\left\{ \left[ \max_{1\leq n\leq N}\left( \left\| q_{n}^{\ast
}\right\| _{C^{1+\alpha }\left( \overline{\Omega }\right) }\right) +2\left\|
T^{\ast }\right\| _{C^{1+\alpha }\left( \overline{\Omega }\right) }+\left\|
\nabla q_{1}^{\ast }\right\| _{C^{\alpha }\left( \overline{\Omega }\right)
}^{2}+1\right] ,C^{\ast },C_{1}\right\} ,
\end{equation*}
where $C^{\ast }$ and $C_{1}$ are constants from (6.9) and (6.10)
respectively. Consider the Dirichlet boundary value problem 
\begin{equation*}
\Delta u+\sum\limits_{j=1}^{3}b_{j}(x)u_{x_{j}}-d(x)u=f\left( x\right) \text{%
, }x\in \Omega ,
\end{equation*}
\begin{equation*}
u\mid _{\partial \Omega }=g\left( x\right) ,g\in C^{2+\alpha }\left(
\partial \Omega \right) ,
\end{equation*}
where functions 
\begin{equation*}
b_{j},d,f\in C^{\alpha }\left( \overline{\Omega }\right) ,d\left( x\right)
\geq 0;\text{ }\max \left( \left\| b_{j}\right\| _{C^{\alpha }\left( 
\overline{\Omega }\right) },\left\| d\right\| _{C^{\alpha }\left( \overline{%
\Omega }\right) }\right) \leq 4M^{\ast }.
\end{equation*}
By the Schauder theorem there exists unique solution $u\in C^{2+\alpha
}\left( \overline{\Omega }\right) $ of this problem and with a constant $%
K=K\left( M^{\ast },\Omega \right) >0$ the following estimate holds 
\begin{equation*}
\left\| u\right\| _{C^{2+\alpha }\left( \overline{\Omega }\right) }\leq K%
\left[ \left\| g\right\| _{C^{2+\alpha }\left( \partial \Omega \right)
}+\left\| f\right\| _{C^{\alpha }\left( \overline{\Omega }\right) }\right] .
\end{equation*}
In Theorem 6.1 we use a function $T_{appr}\left( x,z,\overline{s}\right) $
instead of the above constructed function $T_{1}\left( x,z,\overline{s}%
\right) $ only because the latter was constructed for a rectangle, while
Theorem 6.1 works with a convex bounded domain, also see the beginning of
this section.

\textbf{Theorem 6.1}. \emph{Let }$\Omega \subset \mathbb{R}^{2}$\emph{\ be a
convex bounded domain with the boundary }$\partial \Omega \in C^{3}.$ \emph{%
Suppose that an approximation }$T_{appr}\left( x,z,\overline{s}\right) $ 
\emph{for the tail is constructed in such a way that} 
\begin{equation}
\left\| T_{appr}-T^{\ast }\right\| _{C^{2+\alpha }\left( \overline{\Omega }%
\right) }\leq \xi ,  \tag{6.11}
\end{equation}
\emph{where }$\xi \in \left( 0,1\right) $\emph{\ is a sufficiently small
number and that this function }$T_{appr}\left( x,z,\overline{s}\right) $ 
\emph{is used in (5.3), (5.4) instead of the function} $T_{1}\left( x,z,%
\overline{s}\right) .$ \emph{Denote }$\eta =h+\sigma +\xi +\varepsilon .$ 
\emph{Suppose that the number }$\beta :=\overline{s}-s_{0}=Nh$ \emph{is }$%
such$\emph{\ that } 
\begin{equation}
\beta \leq \frac{1}{48KM^{\ast }}.  \tag{6.12}
\end{equation}
\emph{Then there exists a sufficiently small number }$\eta _{0}=\eta
_{0}\left( K\left( M^{\ast },\Omega \right) ,M^{\ast },c,\underline{s},%
\overline{s}\right) \in \left( 0,1\right) $\emph{\ and a sufficiently large
small number }$h_{0}=h_{0}\left( K\left( M^{\ast },\Omega \right) ,M^{\ast
},c,\underline{s},\overline{s}\right) \in \left( 0,1\right) $ \emph{such
that for all }$\eta \in \left( 0,\eta _{0}\right) $ \emph{and} \emph{for
every integer }$n\in \left[ 1,N\right] $\emph{\ the following estimates hold}
\begin{equation}
\left\| q_{n}-q_{n}^{\ast }\right\| _{C^{2+\alpha }\left( \overline{\Omega }%
\right) }\leq 2KM^{\ast }\left( h+3\eta \right) ,  \tag{6.13}
\end{equation}
\begin{equation}
\left\| q_{n}\right\| _{C^{2+\alpha }\left( \overline{\Omega }\right) }\leq
2M^{\ast }.  \tag{6.14}
\end{equation}
\emph{\ }

\textbf{Remark 6.2.} As it was stated above, unlike the time dependent case
of [3], we cannot prove the estimate (6.11) for tails. However, we observe
convergence of tails in computations if taking $T:=T_{1}$ as in (4.7).

\subsection{Proof of Theorem 6.1}

In the course of this proof we assume that $\eta \in \left( 0,\eta
_{0}\right) ,h\in \left( 0,\eta _{0}\right) .$ Denote 
\begin{equation}
\widetilde{q}_{n,k}(\mathbf{x})=q_{n,k}(\mathbf{x})-q^{\ast }(\mathbf{x}),%
\widetilde{T}(\mathbf{x})=T_{appr}(\mathbf{x})-T^{\ast }(\mathbf{x}),%
\widetilde{\psi }_{n}=\psi _{n}-\psi _{n}^{\ast },  \tag{6.16}
\end{equation}
\begin{equation}
v_{n,k}\left( \mathbf{x},s_{n}\right)
=-hq_{n,k}-h\sum\limits_{j=1}^{n-1}q_{j}+T_{appr}\left( \mathbf{x}\right)
,u_{n,k}\left( \mathbf{x},s_{n}\right) =\exp \left[ v_{n,k}\left( \mathbf{x}%
,s_{n}\right) \right] ,  \tag{6.17}
\end{equation}
\begin{equation}
\widetilde{v}_{n,k}\left( \mathbf{x},s_{n}\right) =v_{n,k}\left( \mathbf{x}%
,s_{n}\right) -v_{n}^{\ast }\left( \mathbf{x},s_{n}\right) ,\widetilde{a}%
_{n,k}\left( \mathbf{x}\right) =a_{nk}\left( \mathbf{x}\right) -a^{\ast
}\left( \mathbf{x}\right) .  \tag{6.18}
\end{equation}
The proof basically consists in estimating these differences.

First, we estimate $\widetilde{q}_{1,1}.$ Set in (6.8) $n=1$ and subtract it
from (5.3) at $k=1$, recalling that $q_{1,0}=0.$ We obtain 
\begin{equation*}
\Delta \widetilde{q}_{1,1}-2\nabla \widetilde{q}_{1,1}\nabla
T_{appr}=2\nabla q_{1}^{\ast }\nabla \widetilde{T}-2h\left( \nabla
q_{1}^{\ast }\right) ^{2}-F_{1}^{\ast },
\end{equation*}
\begin{equation*}
\widetilde{q}_{1,1}\mid _{\partial \Omega }=\widetilde{\psi }_{1}.
\end{equation*}
By Schauder theorem, and (6.9)-(6.11) we obtain 
\begin{equation}
\left\| \widetilde{q}_{1,1}\right\| _{C^{2+\alpha }\left( \overline{\Omega }%
\right) }\leq KM^{\ast }\left( h+\eta \right) .  \tag{6.19}
\end{equation}
Hence, 
\begin{equation}
\left\| \widetilde{q}_{1,1}+q_{1}^{\ast }\right\| _{C^{2+\alpha }\left( 
\overline{\Omega }\right) }=\left\| q_{1,1}\right\| _{C^{2+\alpha }\left( 
\overline{\Omega }\right) }\leq M^{\ast }+KM^{\ast }\left( h+\eta \right)
\leq 2M^{\ast }.  \tag{6.20}
\end{equation}

Now we estimate $\widetilde{q}_{1,k}.$ Set in (6.8) $n=1$ and subtract it
from (5.3)$.$ Then 
\begin{equation}
\Delta \widetilde{q}_{1,k}+2\left( \nabla T_{appr}-h\nabla q_{1,k-1}\right)
\nabla \widetilde{q}_{1,k}=  \tag{6.21}
\end{equation}
\begin{equation*}
-2\nabla q_{1}^{\ast }\left( \nabla \widetilde{T}-2h\nabla \widetilde{q}%
_{1,k-1}\right) -F_{1}^{\ast },
\end{equation*}
and also 
\begin{equation}
\widetilde{q}_{1,k}\mid _{\partial \Omega }=\widetilde{\psi }_{1}. 
\tag{6.22}
\end{equation}
Set in (6.21) $k=2$. Then (6.19) and (6.20) imply that 
\begin{equation}
\left\| 2\left( \nabla T_{appr}-h\nabla q_{1,1}\right) \right\| _{C^{\alpha
}\left( \overline{\Omega }\right) }\leq 2\left( M^{\ast }+2hM^{\ast }\right)
\leq 3M^{\ast },  \tag{6.23}
\end{equation}
\begin{equation*}
\left\| 2\nabla q_{1}^{\ast }\left( \nabla \widetilde{T}-2h\nabla \widetilde{%
q}_{1,1}\right) +F_{1}^{\ast }\right\| _{C^{\alpha }\left( \overline{\Omega }%
\right) }\leq 2M^{\ast }\left[ \xi +2KM^{\ast }h\left( h+\eta \right) +\eta
/2\right] .
\end{equation*}
By (6.12) $2KM^{\ast }h<1/2.$ Hence, 
\begin{equation*}
\left\| 2\nabla q_{1}^{\ast }\left( \nabla \widetilde{T}-2h\nabla \widetilde{%
q}_{1,1}\right) +F_{1}^{\ast }\right\| _{C^{\alpha }\left( \overline{\Omega }%
\right) }\leq 2M^{\ast }\left( h+2\eta \right) .
\end{equation*}
Hence, by (6.21)-(6.23) and Schauder's theorem 
\begin{equation*}
\left\| \widetilde{q}_{1,2}\right\| _{C^{2+\alpha }\left( \overline{\Omega }%
\right) }\leq 2KM^{\ast }\left( h+3\eta \right) .
\end{equation*}
and similarly with (6.19) 
\begin{equation*}
\left\| q_{1,2}\right\| _{C^{2+\alpha }\left( \overline{\Omega }\right)
}\leq 2M^{\ast }.
\end{equation*}

Assume that 
\begin{equation}
\left\| \widetilde{q}_{1,k-1}\right\| _{C^{2+\alpha }\left( \overline{\Omega 
}\right) }\leq 2KM^{\ast }\left( h+3\eta \right) ,\left\| q_{1,k-1}\right\|
_{C^{2+\alpha }\left( \overline{\Omega }\right) }\leq 2M^{\ast }.  \tag{6.24}
\end{equation}
We now estimate the function $\widetilde{q}_{1,k-1}.$ Similarly with the
above 
\begin{equation}
\left\| 2\left( \nabla T_{appr}-h\nabla q_{1,k-1}\right) \right\|
_{C^{\alpha }\left( \overline{\Omega }\right) }\leq 2\left( M^{\ast
}+2hM^{\ast }\right) \leq 6M^{\ast }.  \tag{6.25}
\end{equation}
Next, using (6.24), we obtain 
\begin{equation*}
\left\| 2\nabla q_{1}^{\ast }\left( \nabla \widetilde{T}-2h\nabla \widetilde{%
q}_{1,k-1}\right) +F_{1}^{\ast }\right\| _{C^{\alpha }\left( \overline{%
\Omega }\right) }\leq 2M^{\ast }\left[ \xi +2KM^{\ast }h\left( h+\eta
\right) +\eta /2\right] .
\end{equation*}
By (6.12) $2KM^{\ast }h\left( h+\eta \right) <1/2\left( h+\eta \right) .$
Hence, 
\begin{equation*}
\xi +2KM^{\ast }h\left( h+\eta \right) +\eta /2\leq h+2\eta .
\end{equation*}
Hence, 
\begin{equation}
\left\| 2\nabla q_{1}^{\ast }\left( \nabla \widetilde{T}-2h\nabla \widetilde{%
q}_{1,k-1}\right) +F_{1}^{\ast }\right\| _{C^{\alpha }\left( \overline{%
\Omega }\right) }\leq 2M^{\ast }\left( h+2\eta \right) .  \tag{6.26}
\end{equation}
Hence, Schauder's theorem, (6.21), (6.25) and (6.26) lead to 
\begin{equation}
\left\| \widetilde{q}_{1,k}\right\| _{C^{2+\alpha }\left( \overline{\Omega }%
\right) }\leq 2KM^{\ast }\left( h+3\eta \right) ,\left\| q_{1,k-1}\right\|
_{C^{2+\alpha }\left( \overline{\Omega }\right) }\leq 2M^{\ast },k=1,2,... 
\tag{6.27}
\end{equation}

We now estimate the function $\widetilde{q}_{n,k},$ assuming that (6.27)
holds for functions $\widetilde{q}_{i},q_{i}$ with $j<n,$ as well as for
functions $\widetilde{q}_{n,m},q_{n.m}$ with $m\leq k-1.$ In other words, we
assume that 
\begin{equation}
\left\| \widetilde{q}_{j}\right\| _{C^{2+\alpha }\left( \overline{\Omega }%
\right) }\leq 2KM^{\ast }\left( h+3\eta \right) ,\left\| q_{j}\right\|
_{C^{2+\alpha }\left( \overline{\Omega }\right) }\leq 2M^{\ast }.  \tag{6.28}
\end{equation}
and 
\begin{equation}
\left\| \widetilde{q}_{n,m}\right\| _{C^{2+\alpha }\left( \overline{\Omega }%
\right) }\leq 2KM^{\ast }\left( h+3\eta \right) ,\left\| q_{n,m}\right\|
_{C^{2+\alpha }\left( \overline{\Omega }\right) }\leq 2M^{\ast },m\leq k-1. 
\tag{6.29}
\end{equation}
Subtracting (6.8) from (5.4) and using (5.5), (6.3) and (6.16), we obtain 
\begin{equation}
\Delta \widetilde{q}_{n.k}-2\left( h\sum\limits_{j=1}^{n-1}\nabla
q_{j}-\left( \nabla T_{appr}-2h\nabla q_{n,k-1}\right) \right) \nabla 
\widetilde{q}_{n.k}=  \tag{6.30}
\end{equation}
\begin{equation*}
2\left( h\sum\limits_{j=1}^{n-1}\nabla \widetilde{q}_{j}\right) \nabla
q_{n}^{\ast }+2\nabla q_{n}^{\ast }\nabla \widetilde{T}+2h\nabla q_{n}^{\ast
}\nabla \widetilde{q}_{n.k-1},
\end{equation*}
\begin{equation}
\widetilde{q}_{n.k}\mid _{\partial \Omega }=\widetilde{\psi }_{n}. 
\tag{6.31}
\end{equation}
Estimate first the coefficient at $\nabla \widetilde{q}_{n.k}$ in (6.30).
Using (6.28) and (6.29), we obtain 
\begin{equation*}
2\left\| h\sum\limits_{j=1}^{n-1}\nabla q_{j}\right\| _{C^{\alpha }\left( 
\overline{\Omega }\right) }\leq 4M^{\ast }Nh,
\end{equation*}
\begin{equation*}
2\left\| \nabla T_{appr}-2h\nabla q_{n,k-1}\right\| _{C^{\alpha }\left( 
\overline{\Omega }\right) }\leq 2\left( M^{\ast }+2hM^{\ast }\right) \leq
3M^{\ast }.
\end{equation*}
Since by (6.12) $4M^{\ast }\overline{N}h\leq M^{\ast },$ then the estimate
for that coefficient is 
\begin{equation}
2\left\| h\sum\limits_{j=1}^{n-1}\nabla q_{j}\right\| _{C^{\alpha }\left( 
\overline{\Omega }\right) }+2\left\| \nabla T_{appr}-2h\nabla
q_{n,k-1}\right\| _{C^{\alpha }\left( \overline{\Omega }\right) }\leq
4M^{\ast }.  \tag{6.32}
\end{equation}
Hence, we can apply Schauder's theorem with the constant $K.$ Now we
estimate the right hand side of equation (6.30). Using (6.28) and (6.29), we
obtain 
\begin{equation*}
\left\| 2\left( h\sum\limits_{j=1}^{n-1}\nabla \widetilde{q}_{j}\right)
\nabla q_{n}^{\ast }\right\| _{C^{\alpha }\left( \overline{\Omega }\right)
}\leq 4K\left( M^{\ast }\right) ^{2}Nh\left( h+3\eta \right) ,
\end{equation*}
\begin{equation*}
\left\| 2h\nabla q_{n}^{\ast }\nabla \widetilde{q}_{n.k-1}\right\|
_{C^{\alpha }\left( \overline{\Omega }\right) }\leq 4K\left( M^{\ast
}\right) ^{2}h\left( h+3\eta \right) ,
\end{equation*}
\begin{equation*}
\left\| 2\nabla q_{n}^{\ast }\nabla \widetilde{T}\right\| _{C^{\alpha
}\left( \overline{\Omega }\right) }\leq 2M^{\ast }\eta .
\end{equation*}
Estimate the right hand sides of the last three inequalities. By (6.12) we
have 
\begin{equation*}
4K\left( M^{\ast }\right) ^{2}\overline{N}h\left( h+3\eta \right) \leq \frac{%
M^{\ast }}{12}\left( h+3\eta \right) ,
\end{equation*}
\begin{equation*}
4K\left( M^{\ast }\right) ^{2}h\left( h+3\eta \right) \leq \frac{M^{\ast }}{%
12}\left( h+3\eta \right) .
\end{equation*}
Hence, 
\begin{equation}
\left\| 2\left( h\sum\limits_{j=1}^{n-1}\nabla \widetilde{q}_{j}\right)
\nabla q_{n}^{\ast }\right\| _{C^{\alpha }\left( \overline{\Omega }\right)
}+\left\| 2h\nabla q_{n}^{\ast }\nabla \widetilde{q}_{n.k-1}\right\|
_{C^{\alpha }\left( \overline{\Omega }\right) }+\left\| 2h\nabla q_{n}^{\ast
}\nabla \widetilde{q}_{n.k-1}\right\| _{C^{\alpha }\left( \overline{\Omega }%
\right) }  \tag{6.33}
\end{equation}
\begin{equation*}
\leq \frac{M^{\ast }}{6}\left( h+3\eta \right) +2M^{\ast }\eta =M^{\ast
}\left( \frac{h}{6}+\frac{3}{2}\eta \right) .
\end{equation*}
By (6.30)-(6.33) and Schauder theorem we obtain 
\begin{equation}
\left\| \widetilde{q}_{n.k}\right\| _{C^{2+\alpha }\left( \overline{\Omega }%
\right) }\leq KM^{\ast }\left( \frac{h}{6}+\frac{3}{2}\eta \right) +K\eta
\leq KM^{\ast }\left( \frac{h}{6}+\frac{5}{2}\eta \right) \leq KM^{\ast
}\left( h+3\eta \right) .  \tag{6.34}
\end{equation}
Hence, 
\begin{equation}
\left\| q_{n.k}\right\| _{C^{2+\alpha }\left( \overline{\Omega }\right)
}=\left\| \widetilde{q}_{n.k}+q_{n}^{\ast }\right\| _{C^{2+\alpha }\left( 
\overline{\Omega }\right) }\leq KM^{\ast }\left( h+3\eta \right) +M^{\ast
}\leq 2M^{\ast }.  \tag{6.35}
\end{equation}
Estimates (6.34) and (6.35) complete the proof of this theorem. $\square $

\section{Numerical Studies}

We have performed numerical experiments on several cases of reconstructions
using the method discussed above. We have chosen the range of geometrical
parameters of the rectangle $\Omega $, which is typical for optical imaging
of small animals and have chosen the range of optical parameters typical for
biological tissues [1],[11],[13].

\subsection{\protect\bigskip Some Details of Numerical Studies}

For the forward problem, we calculate the solution of the diffusion equation 
\begin{equation}
D\Delta u-\mu _{a}(x,z)u=-\delta \left( x-s,z-z_{m}\right)  \tag{7.1}
\end{equation}
with the conventional condition at the infinity 
\begin{equation}
\underset{\left| \left( x,z\right) \right| \rightarrow \infty }{\lim }%
u(x,z,s)=0,  \tag{7.2}
\end{equation}
where $D=1/\left( 3\mu _{s}^{\prime }\right) \equiv const.>0$ is the
diffusion coefficient, where optical coefficients $\mu _{s}^{\prime }$ and $%
\mu _{a}(x,z)$ were discussed in Introduction (Section 1). In our
computations the function $\mu _{a}(x,z)$ is unknown and \ the constant $\mu
_{s}^{\prime }$ is given. Thus, in our case $a(x,z)=3\mu _{s}^{\prime }\cdot
\mu _{a}(x,z)$ (compare with (1.4)). Consider the rectangle $\Omega ,$%
\begin{equation*}
\Omega =\left\{ (x,z):5cm<x<15cm,5cm<z<10cm\right\} .
\end{equation*}
We assume that 
\begin{equation}
a(x,z)=k^{2}=const.>0\text{ in }\mathbb{R}^{2}\diagdown \Omega .  \tag{7.3}
\end{equation}
We assume that in (7.1) the source position $\left( s,z_{m}\right) $ is
running along the right side of $\Omega ,$ i.e., $z_{m}=L=10cm.$ Also,
consider a bigger rectangle 
\begin{equation*}
\Omega _{0}=\{(x,z):0cm<x<20cm,0cm<z<15cm\}.
\end{equation*}
The reason why we consider the rectangle $\Omega _{0}$ along with the
rectangle $\Omega $ is that it is natural to approximate the solution of the
problem (7.1), (7.2) in the infinite domain by the solution of equation
(7.1) in $\Omega _{0}$ with Robin boundary conditions at $\partial \Omega
_{0}.$ We have established numerically that for the range of parameters we
use the solution of the problem (7.1), (7.2) is close in $\Omega $ to the
solution of equation (7.1) in the bigger rectangle $\Omega _{0}$ with the
Robin boundary conditions at its sides.\ Figure 3 illustrates rectangles $%
\Omega _{0}$ and $\Omega .$

The light sources are located in several positions $\left( x_{i},z\right)
=\left( s_{i},10\right) $ along the right-hand side of the smaller rectangle 
$\Omega $, and receivers, which mimic the so-called CCD camera are located
at the left-hand side of $\Omega $. CCD stands for a ``charge-coupled
device''. A CCD camera is an image sensor, consisting of an integrated
circuit containing an array of linked, or coupled, light-sensitive
capacitors. A typical CCD camera can take up to $512\times 512$ data points
simultaneously, which will provide an adequate amount of data for our
reconstruction. In all three examples, we have used an ideal light source
modeled by the function $-\delta (x-s_{i},z-10)$ in the 2D case of (1.1). In
numerical simulation $\delta (x-s_{i},z-10)=c\eta (s_{i},10)$, where $\eta $
is the finite element at the location, and $c$ is the scaling constant to
ensure that the area equals one.

We use three (3) sources to construct an approximation of the tail functions
which was described above. Next, we use all five (5) sources for the above
layer stripping procedure both in the $s$-derivative and the $s$-integral.

We have generated the data for the forward problem for total of five (5)
different locations of the light source, $s_{i}=1,...,5,$ where $%
s_{1}=0,s_{i}=s_{i-1}+0.625cm,i=2,...,5.$ Hence, $N=4$ and we have used four
(4) functions $q_{n}.$ \ An increase of the number $N$ \ did not result in
significant improvements of results. Note that a similar observation took
place in\ numerical experiments of [3] for the case of an increase of the
number of functions $q_{n}$ \ after a certain ``limit''. In our
reconstruction method, we use the solution of the forward problem to
generate the data for the inverse, add noise to the measurement data, and
reconstruct the absorption coefficient $\mu _{a}(x,z)$ in $\Omega $. The
domain $\Omega $ will be our basic computational domain for our inverse
calculations. 

In our examples, the coefficients in equation (7.1) are $D=0.02cm$ uniformly
and $\mu _{a}=0.1cm^{-1}$ at all grids except off the inclusions, and in
inclusions $\mu _{a}$ ranges from 0.1 to 0.3 $cm^{-1}$. The maximum
inclusion/background contrast is 3:1 in our computations. Our algorithm
calculates the forward problem with Robin boundary conditions at $\partial
\Omega _{0}$, given the distribution of the absorption coefficient. A total
of 130$\times $93 rectangular finite elements is used for forward
calculations for the domain $\Omega _{0}$.

For the simulated boundary measurements we take the solution of the forward
problems along the left and lower boundaries of $\Omega $ to construct first
guess for tails (section 4.1) and then we take all 4 sides for the rest of
the problem. The number of measuring points is 65 along the left edge of $%
\Omega $ and 31 along the lower edge of $\Omega $. The number of measuring
points at the low left corner is shared by both sides and therefore the
total number of independent measuring points is 95.

For each detector position, we introduce the random noise as the random
process with respect to the detector locations, $\widetilde{\varphi }(%
\mathbf{x},s_{k})=\varphi (\mathbf{x},s_{k})\left[ 1+\chi \left( \mathbf{x}
\right) \right] ,$ where $\chi \left( \mathbf{x}\right) $ is the random
variable, which we introduce as $\chi =0.02W$, where $W$ is a white noise
with equal distribution at [-1,1]. Hence, this is 2\% of the multiplicative
random noise. To obtain the realistic first $s$-derivatives, we started with
simulated light distribution $u(x,z,s)$ added with similar noise to simulate
the situation used in applications and let $v=\ln u$. Then we take first
derivatives with respect to $s$ as shown in the paragraph below.

A regularization method was introduced to pre-process the noise in the
measurement data. We use a polynomial approximation with respect to the
detector location $\mathbf{x}$. In our setting, the measurements are
collected at 65 locations along the left boundary and 31 points along the
lower boundary. We use an eight order polynomial to approximate functions $%
\widetilde{\varphi }(\mathbf{x},s_{k})$ with respect to $\mathbf{x}$\textbf{%
\ }for each $s_{k}.$ The polynomial is optimal in the least square sense
[12],[14], and its sub-routing is commonly available, see for example

\text{\textit{http://perso.orange.fr/jean-pierre.moreau/f\_lstsqr.html}}. We
demonstrate the essence of the approximation in Figure 4. Thus, we have
obtained approximate polynomial functions $\overline{\varphi }(x,s_{k}).$ We
use functions $\overline{\varphi }(x,s_{k})$ instead of $\varphi (x,s_{k}).$
The first $s$-derivatives are processed afterwards by the formula 
\begin{equation*}
f^{\prime }(s_{1})\approx \frac{f\left( s_{2}\right) -f\left( s_{1}\right) }{%
s_{2}-s_{1}}.
\end{equation*}
and similar ones for other source locations.

\subsection{Numerical Experiments}

In the following numerical examples, we illustrate the results in a few
different shapes and locations of the two inclusions. Our method has shown
its success in dealing with those cases. In all cases, iterations with
respect to tail (as described in section 4.2) were only done for the first
light sources. Additional light sources did not bring any significant
changes to reconstruction, therefore iterations with respect to tails are
not shown in this paper. The total number of elements K in (3.8) is 450 
in our calculation.

The Convergence Criterion for functions $a_{m}$ in the procedure of finding
the second accelerator $T_{1}$ for tails is 
\begin{equation}
||a_{1m_{1}}(x,z)-a_{1m_{1}-1}(x,z)||\equiv {\frac{\sqrt{\sum%
\limits_{i=1,...,i_{max},j=1,...,j_{max}}|(a_{1,m_{1}}(x_{i},z_{j})-a_{1,m_{1}-1}(x_{i},z_{j}))|^{2}%
}}{\sqrt{N_{1}}\max {|a_{1,m_{1}-1}(x_{i},z_{j})|}}}\leq \varepsilon , 
\tag{7.1}
\end{equation}
where $N_{1}=i_{max}j_{max}$ is the total number of finite elements. In all
our examples, $\varepsilon =10^{-5}$. The number of iterations required for
convergence is listed below:

\textit{\ The number of iterations required for convergence}

\begin{tabular}{|r|r|r|}
\hline
Example 1 & Example 2 & Example 3 \\ \hline
52 & 50 & 48 \\ \hline
\end{tabular}

\textbf{Example 1}. Inclusions are two circles with the radius 1 cm, and
their centers are placed 2 cm off the left edge. The coefficient is $\mu
_{a}(x,z)=0.3$ inside inclusion and $\mu _{a}(x,z)=0.1=k^{2}$ outside of
inclusions. We have also added 2\% of random noise to the boundary
measurements, see subsection 7.1.

Figure 5a displays the original distribution and its 1-d cross section.
Figure 5b shows reconstruction from the noisy data and its 1-d cross section.

The relative errors of the reconstruction are as follows:

\textit{Table 1. The relative errors of reconstructions in Example 1}

\begin{tabular}{|r|r|r|}
\hline
$RMSE$ & $AME$ & $ME$ \\ \hline
0.312366805537619 & 0.115817263155519 & -0.058254402556204 \\ \hline
\end{tabular}

Note that for the data set $(x_{1},x_{2},\cdots ,x_{N_{1}})$ and its
approximation $(\hat{x}_{1},\hat{x}_{2},\cdots ,\hat{x}_{N_{1}})$, the
values of the function $\mu _{a}(x,z)$ taken at each of the grid points, the
Relative Root Mean Square Error (RMSE), Relative Absolute Mean Error (MAE)
and Relative Mean Error (ME) are calculated by 
\begin{equation*}
RMSE={\frac{\sqrt{\sum\limits_{k=1}^{N_{1}}(x_{k}-\hat{x}_{k})^{2}}}{\sqrt{%
N_{1}}\max {|x_{k}|}}},MAE={\frac{\sum\limits_{k=1}^{N_{1}}|x_{k}-\hat{x}%
_{k}|}{N_{1}\max {|x_{k}|}}},ME={\frac{\sum\limits_{k=1}^{N_{1}}(x_{k}-\hat{x%
}_{k})}{N_{1}\max {|x_{k}|}}}.
\end{equation*}

In our case $x_{k}$ are correct values of the coefficient $\mu _{a}\left(
x,z\right) $ at the grid points of the sub-rectangle 
\begin{equation*}
\Omega ^{\prime }=\left\{ \left( x,z\right) :5cm<x<15cm,5cm<z<8cm\right\}
\subset \Omega .
\end{equation*}

We illustrate in Figure 5c the difference of two consecutive reconstruction 
\begin{equation*}
||a_{m}(x,z)-a_{m-1}(x,z)||\equiv {\frac{\sqrt{\sum%
\limits_{i=1,...,i_{max},j=1,...,j_{max}}|(a_{m}(x_{i},z_{j})-a_{m-1}(x_{i},z_{j}))|^{2}%
}}{\sqrt{N_{1}}\max {|a_{m-1}(x_{i},z_{j})|}}}
\end{equation*}
as a function of the number of iteration $m$. Figure 5d depicts the relative
error in comparison with actual inclusion expressed by 
\begin{equation*}
RMSE\equiv {\frac{\sqrt{\sum%
\limits_{i=1,...,i_{max},j=1,...,j_{max}}|(a_{m}(x_{i},z_{j})-a(x_{i},z_{j}))|^{2}%
}}{\sqrt{N_{1}}\max {|a(x_{i},z_{j})|}}}
\end{equation*}
as a function of the number of iterations $m$.

\textbf{Example 2}. Inclusions are two circles of the radius 1cm, and their
centers are placed 2 cm off the left edge. The coefficient $\mu _{a}(x,z)$
is defined as 
\begin{equation}
\mu _{a}(x,z)=\left\{ 
\begin{array}{c}
\max \left[ 0.3\cos d(x,z),0.1\right] ,\text{ inside of each circle} \\ 
0.1\text{ otherwise,}
\end{array}
\right\}  \tag{7.4}
\end{equation}
where $d(x,z)$ is the minimum distance to center of each of these two
circles.

Figure 6a displays the original function in \ two inclusions and its 1-d
cross section. 
Figure 6b shows the reconstruction result with 2\% noise and its 1-d cross
section.

The relative errors of reconstruction are as follows:

\textit{Table 2. The relative errors of reconstruction in Example 2}

\begin{tabular}{|r|r|r|}
\hline
$RMSE$ & $AME$ & $ME$ \\ \hline
0.261946376827213 & 0.087514613280764 & -0.024755081675791 \\ \hline
\end{tabular}

We illustrate in Figure 6c the difference of two consecutive reconstruction $%
||a_{m}(x,z)-a_{m-1}(x,z)||$ as a function of the number of iteration $m$.
Figure 6d depicts the relative error in comparison with actual inclusion
expressed by $RMSE$ as a function of the number of iterations $m$.

\textbf{Example 3}. Inclusions are two circles of the radius 1 cm and 0.6
cm, whose centers are placed 2 cm off the left edge. The coefficient $\mu
_{a}(x,z)$ is defined as 
\begin{equation}
\mu _{a}(x,z)=\left\{ 
\begin{array}{c}
\max \left[ 0.3(\cos d(x,z)(1+0.1\eta (x,z)),0.1\right] ,\text{ inside of
each circle} \\ 
0.1\text{ otherwise.}
\end{array}
\right\}   \tag{7.5}
\end{equation}
Similarly with (7.4) $d(x,z)$ is the distance to center of the circle. In
(7.5) $\eta $ is a realization of a white noise valued between [-1, 1]. The
random pattern is introduced to test the ability of our method to handle
complex shapes. See Figures 7a,b for results.

\textit{Table 3. The relative errors of reconstruction in Example 3}

\begin{tabular}{|r|r|r|}
\hline
$RMSE$ & $AME$ & $ME$ \\ \hline
0.336142461513025 & 0.084080392591033 & -0.026699370402152 \\ \hline
\end{tabular}

We illustrate in Figure 7c the difference of two consecutive reconstruction $%
||a_{m}(x,z)-a_{m-1}(x,z)||$ as a function of the number of iteration $m$.
Figure 7d depicts the relative error in comparison with actual inclusion
expressed by $RMSE$ as a function of the number of iterations $m$.

\begin{center}
\textbf{Acknowledgments}
\end{center}

The work of all authors was supported by the National Institutes of Health
grant \#1R21Ns052850-01A1. The work of MK was also supported by the U.S.
Army Research Laboratory and U.S. Army Research Office under contract/ grant
number W911NF-05-1-0378.

\begin{center}
\textbf{References}
\end{center}

1. S. Arridge, Optical tomography in medical imaging, Inverse Problems, 15,
841-893, 1999.

2. A.B. Bakushinsky, T.\ Khan and A. Smirnova, Inverse problem in optical
tomography and its numerical investigation by iteratively regularized
methods, J. Inverse and Ill-Posed Problems, 13, 537-551, 2005.

3. L. Beilina and M.V. Klibanov, A globally convergent numerical method for
some coefficient inverse problems with resulting second order elliptic
equations, submitted for publication, a preprint is available on line at
http://www.ma.utexas.edu/mp\_arc/index-07 (preprint number 07-311), at
http://www.math.ntnu.no/preprint/numerics/2007,

and at http://www.math.uncc.edu/people/research/mklibanv.php

4. Yu.A.\ Grazin, M.V.\ Klibanov and T.R.\ Lucas, Numerical solution of a
subsurface imaging inverse problem, \emph{SIAM J. Appl. Math}., 62, 664-683,
2001.

5. E. \ Haber, U.M. \ Asher and D.\ Oldenburg, On optimization techniques for
solving nonlinear inverse problems, Inverse Problems, 16, 1263-1280, 2000.

6. M.V. \ Klibanov and A. \ Timonov, \emph{Carleman Estimates for Coefficient
Inverse Problems and Numerical Applications}, VSP, Utrecht, 2004.

7. O.A. \ Ladyzhenskaya and N.N. \ Uralceva, \emph{Linear and Quasilinear
Elliptic Equations}, Academic Press,\ New York, 1969.

8. H.\ Shan, M.V. \ Klibanov, H.\ Liu, N. \ Pantong and J. \ Su, Numerical
implementation of the convexification algorithm for an optical diffusion
tomograph, \emph{Inverse Problems}, 24, paper number 025026, 2008.

9. J. \ Su, H.\ Shan, H.\ Liu and M.V.\ Klibanov, Reconstruction method from a
multiple-site continuous-wave source for three-dimensional optical
tomography, \emph{J.\ Optical Society of America A}, 23, 2388-2395, 2006.

10. A.N.\ Tikhonov and V.Ya.\ Arsenin, \emph{Solutions of Ill-Posed Problems},
Winston\& Sons, Washington, DC 1977.

11. R.R.\ Alfano, R.R.\ Pradhan  and G.C.\ Tang,  Optical spectroscopic
diagnosis of cancer and normal breast tissues, \emph{J. Opt. Soc. Am. B}, 
6, 1015-1023, 1989.

12. E.W.\ Cheney, \emph{Introduction to Approximation Theory}, 
Chelsea Pub. Co,\ New York, 1982.

13. D. \ Grosenick, H.\ Wabnitz, H.H.\ Rinneberg, K.T.\ Moesta and P.M.\ Schlag 
Development of a time-domain optical mammograph and first \emph{in vivo}
applications, \emph{Applied Optics}, 38, 2827-2943, 1999

14. E.\ Isaacson and H.B.\ Keller, \emph{Analysis of Numerical Methods,}
\ Wiley,\ New York, 1966.

15. R. G.\ Novikov, Multidimensional inverse spectral problem for the equation $-\delta \psi + (v(x) - Eu(x))\psi = 0$,
\emph{ Funk. Anal. Pril.} 22, 11–22 (in Russian), 1988.

16.   R. G.\ Novikov, The inverse scattering problem on a fixed energy level for the two-dimensional Schrodinger
operator, \emph{J. Funk. Anal.}, 103, 409–63, 1992.
 
17. A.\ Nachman,   Global uniqueness for a two-dimensional inverse boundary value problem, \emph{ Ann. Math.},
143, 71–96,1996.

18. P.G.\ Grinevich,   The scattering transform for the two-dimensional operator with a potential that decreases
at infinity at fixed nonzero energy, \emph{ Russ. Math. Surv.}, 55, 3–70,2000.

19.  M.I. \ Belishev,    Boundary control in reconstruction of manifolds and metrics (the BC method),  \emph{ Inverse
Problems},  13, R1–45, 1997.

20. S.I.\ Kabanikhin, A.D.\ Satybaev and M.A. Shishlenin,  \emph{ Direct Methods of Solving Multidimensional Inverse
Hyperbolic Problems}, \ Utrecht: VSP,  2004.

21. V.A.\ Burov, S.A.\ Morozov and O.D.\ Rumyantseva,  Reconstruction of fine-scale structure of acoustical
scatterers on large-scale contrast background, \emph{  Acoust. Imaging},  26, 231–8, 2002.

22. S.\ Siltanen, J.L.\ Mueller and D.\ Isaacson,  An implementation of the reconstruction algorithm of A. Nachman
for the 2-D inverse conductivity problem, \emph{ Inverse Problems},  16, 681–9, 2000.

23. M.I.\ Belishev and V. Yu.\ Gotlib,  Dynamical variant of the BC-method: theory and numerical testing, 
\emph{ J. Inverse Ill-Posed Problems}, 7,  221–40, 1999.

24. M.V.\ Klibanov, A.\ Timonov, Numerical studies on the globally convergent
convexification algorithm in 2D, \emph{ Inverse Problems},  23, 123-138, 2007.

\newpage 

\begin{figure}[!h]
\begin{picture}(0,0)\end{picture}
\includegraphics[width=2in,height=2in,angle=0]{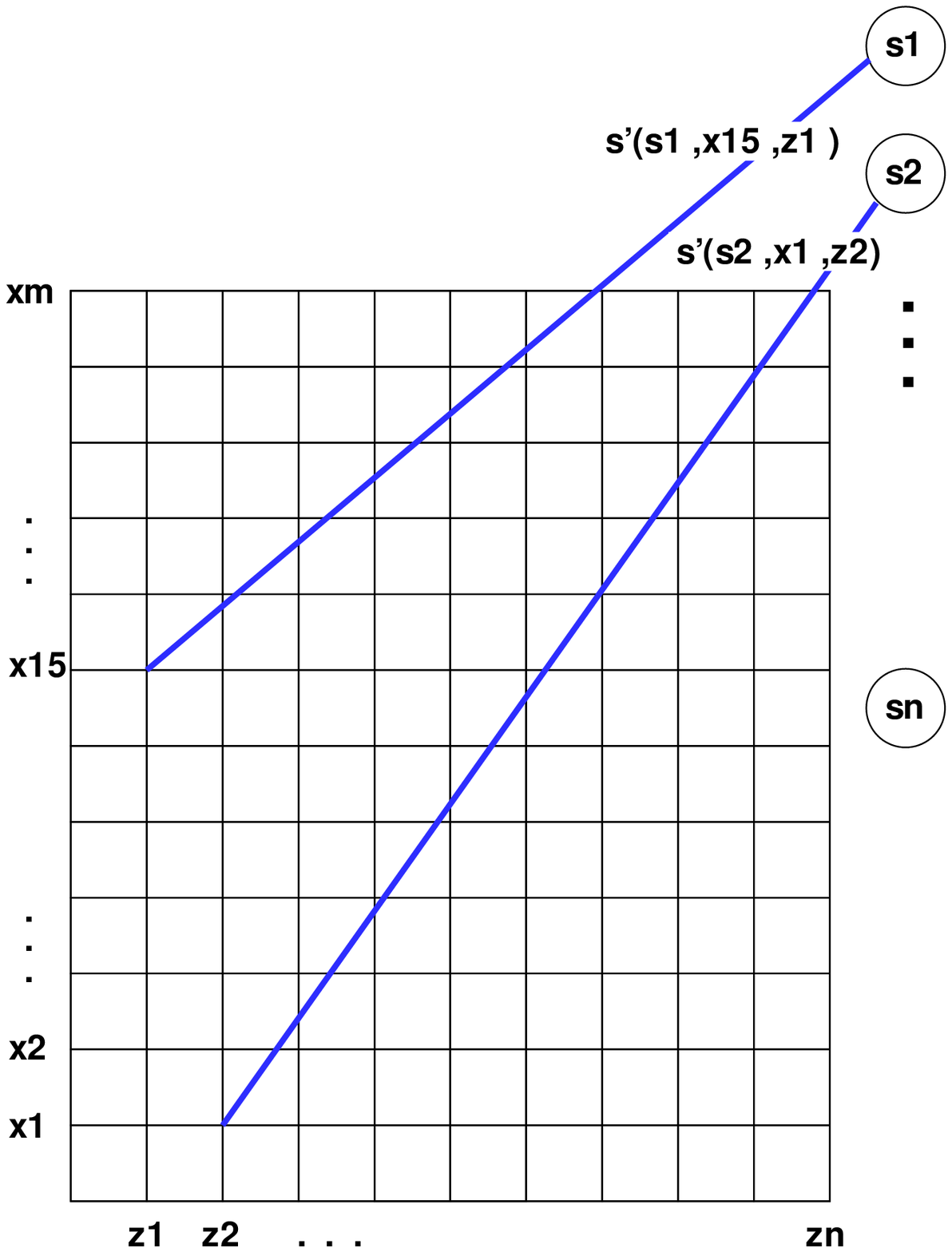} %
\setlength{\unitlength}{1in} 
\hspace{0.5in} \begin{picture}(0,0)\end{picture}
\end{figure}

\textit{Figure 1. The figure illustrates construction of the real distance $%
S $ in the asymptotic expansion of the tail function. The distance $S$
measures the distance to any interior point individually, rather than the
distance to the edge denoted by $s$.} \vskip2in

\begin{figure}[!h]
\begin{picture}(0,0)\end{picture}
\includegraphics[width=2in,height=2in,angle=0]{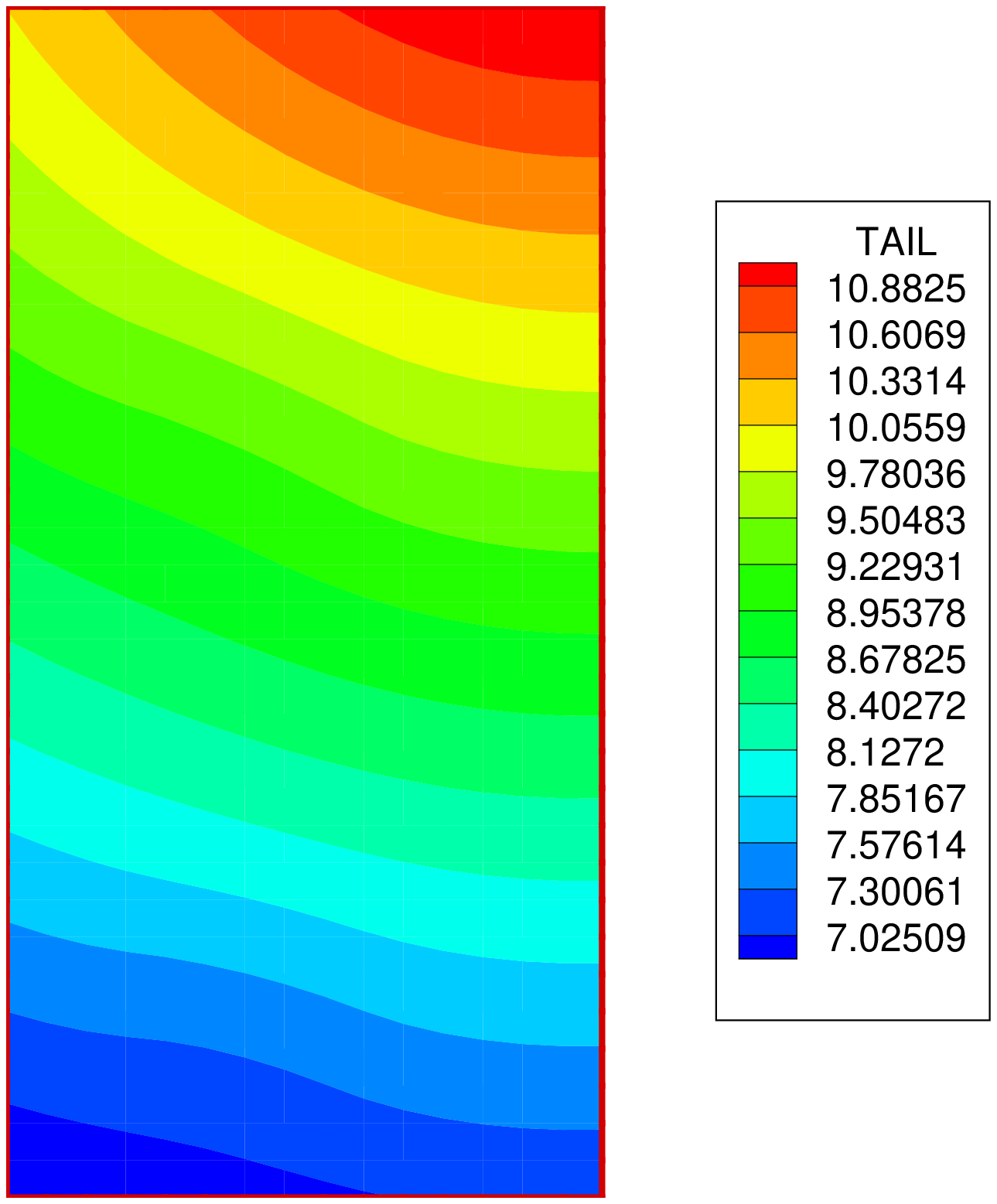} %
\setlength{\unitlength}{1in} 
\hspace{0.5in} \begin{picture}(0,0)\end{picture}
\includegraphics[width=2in,height=2in,angle=0]{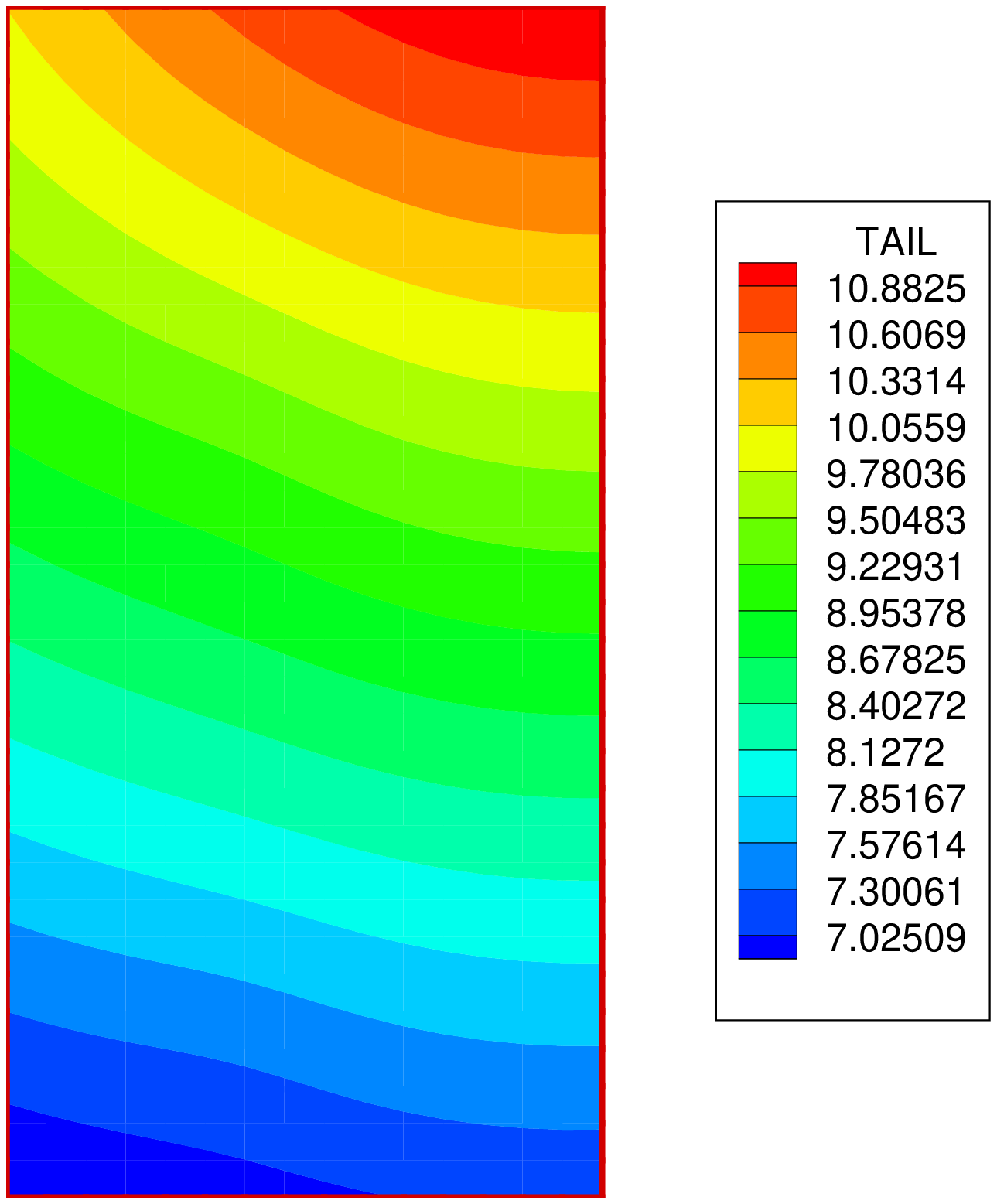} 
\end{figure}

\textit{Figure 2. These two figures illustrate the comparison of actual tail
function $T\left( x,z,\overline{s}\right) $ from forward problem (left
panel) and the calculated tail function $T_{1}\left( x,z,\overline{s}\right) 
$ derived from iterative procedure} (right panel).

\newpage

\begin{figure}[!h]
\begin{picture}(0,0)\end{picture}
\includegraphics[width=2in,height=2in,angle=0]{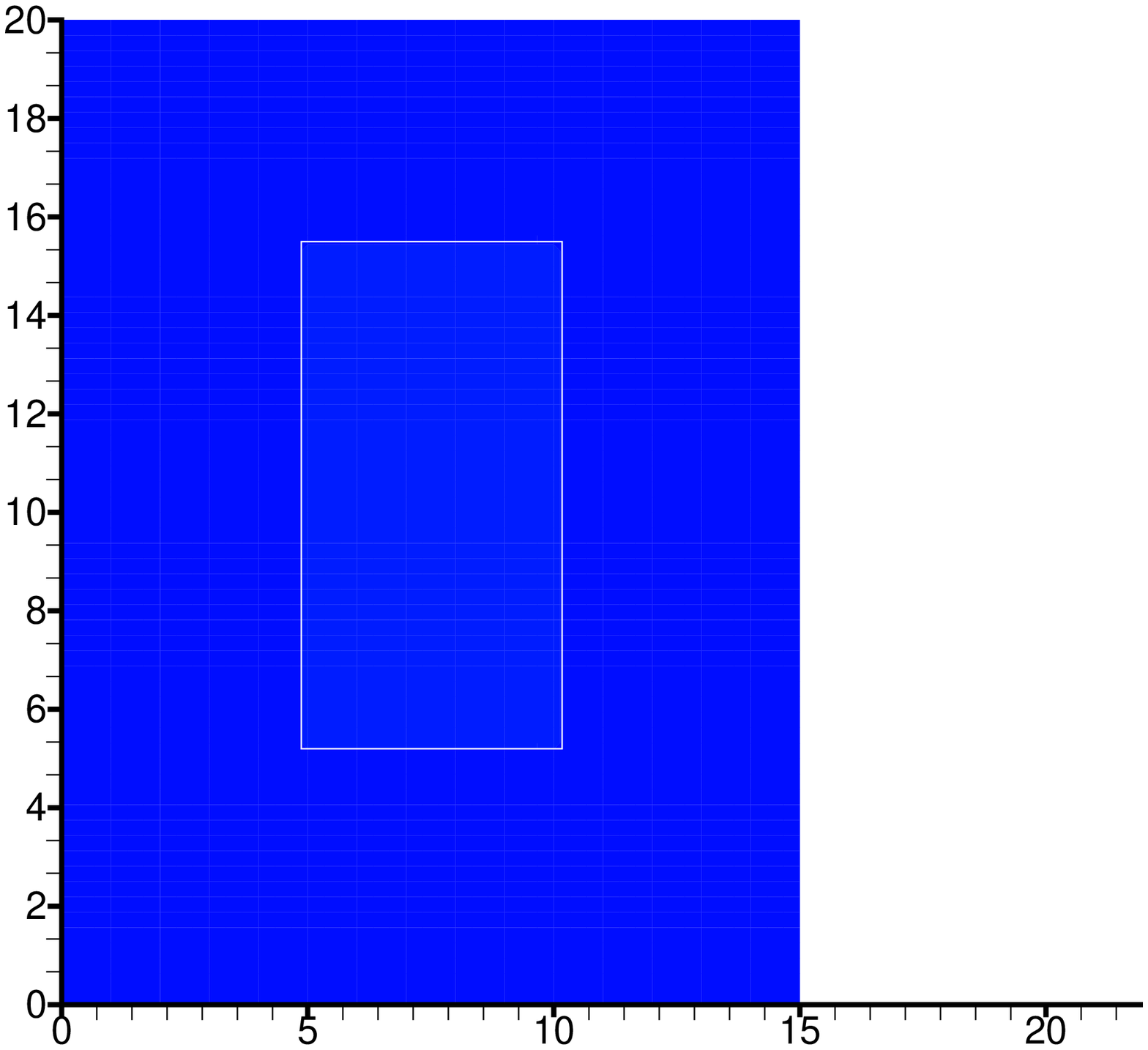} %
\setlength{\unitlength}{1in} 
\hspace{0.5in} \begin{picture}(0,0)\end{picture}
\end{figure}

\textit{Figure 3. The figure illustrates relative positions of rectangles $%
\Omega _{0}$ and $\Omega $. The light source location is at the right edge $%
\Omega $. }


\begin{figure}[!h]
\begin{picture}(0,0)\end{picture}
\includegraphics[width=2in,height=2in,angle=0]{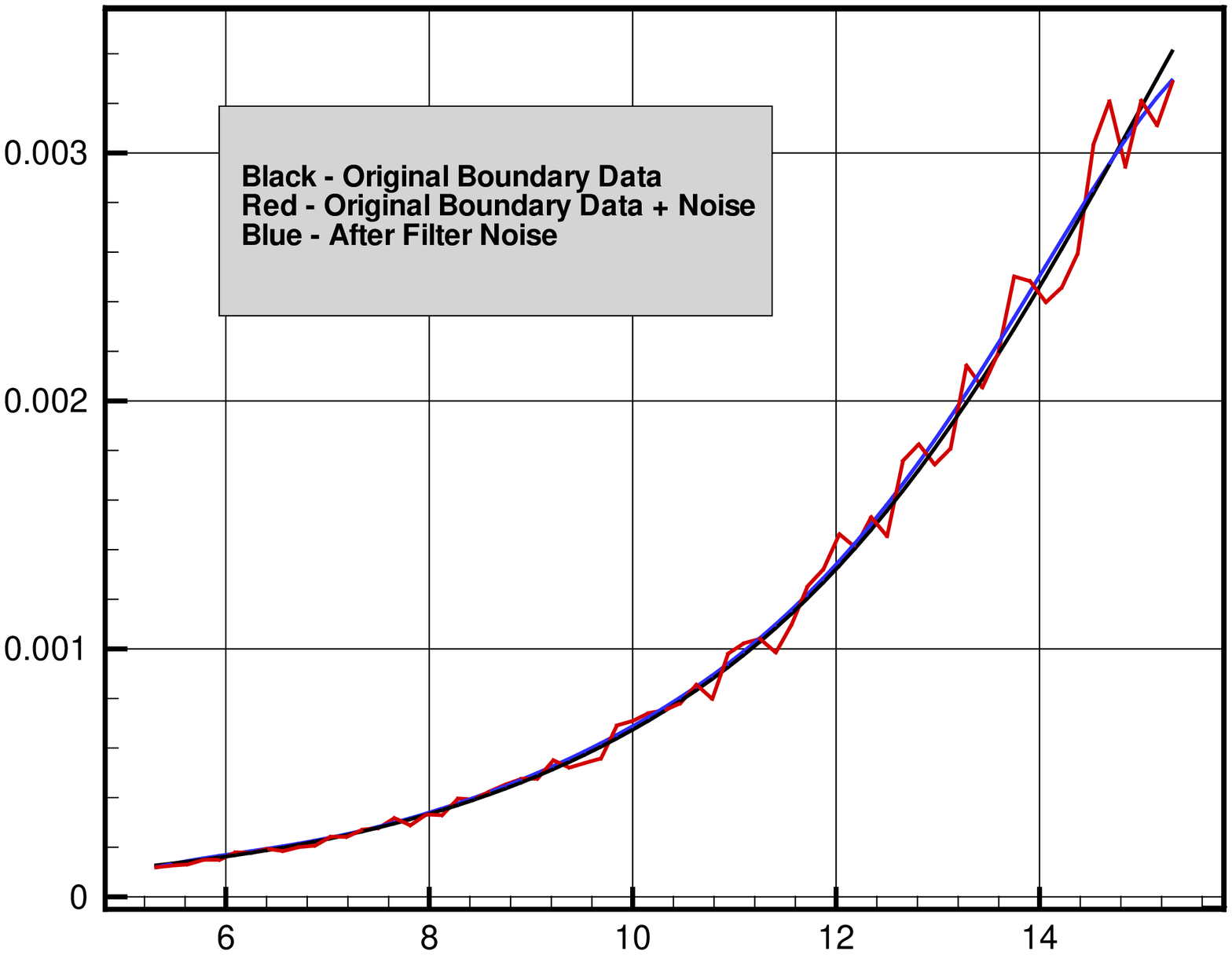} %
\setlength{\unitlength}{1in} 
\hspace{0.5in} \begin{picture}(0,0)\end{picture}
\end{figure}

\textit{\ Figure 4. We use an eighth order polynomial to approximate
functions $\widetilde{\varphi }(x,s_{k})$ and $\widetilde{\psi }(x,s_{k})$
with respect to $x$ for each $s_{k}.$ The polynomial is optimal in the least
square sense, and its sub-routing is commonly available [ 11,12 ]. We
demonstrate the essence of the approximation in Figure 4.}

\newpage 

\begin{figure}[!h]
\begin{picture}(0,0)\end{picture}
\includegraphics[width=2in,height=2in,angle=0]{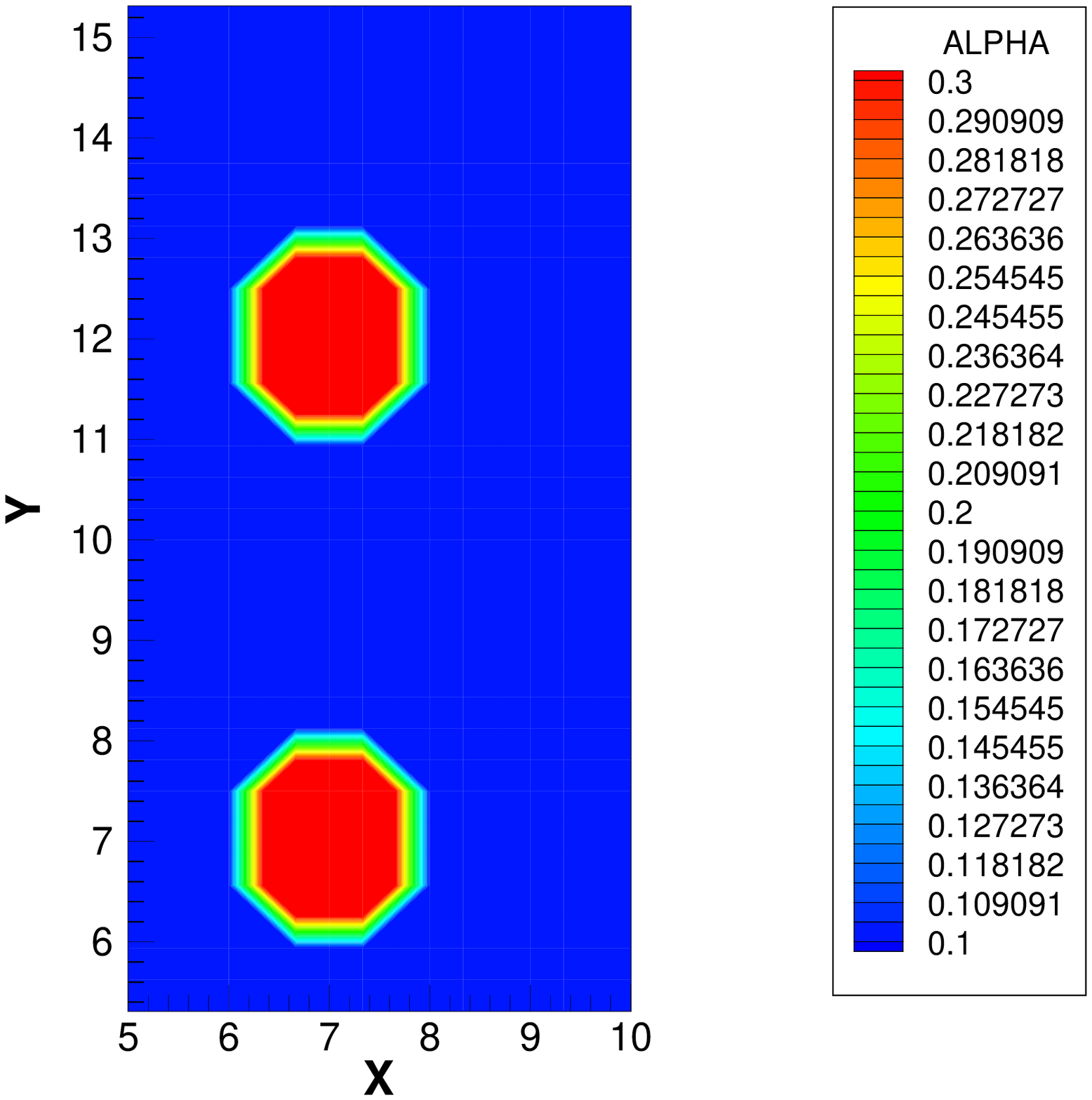} \setlength{%
\unitlength}{1in} 
\hspace{0.5in} \begin{picture}(0,0)\end{picture}
\includegraphics[width=2in,height=2in,angle=0]{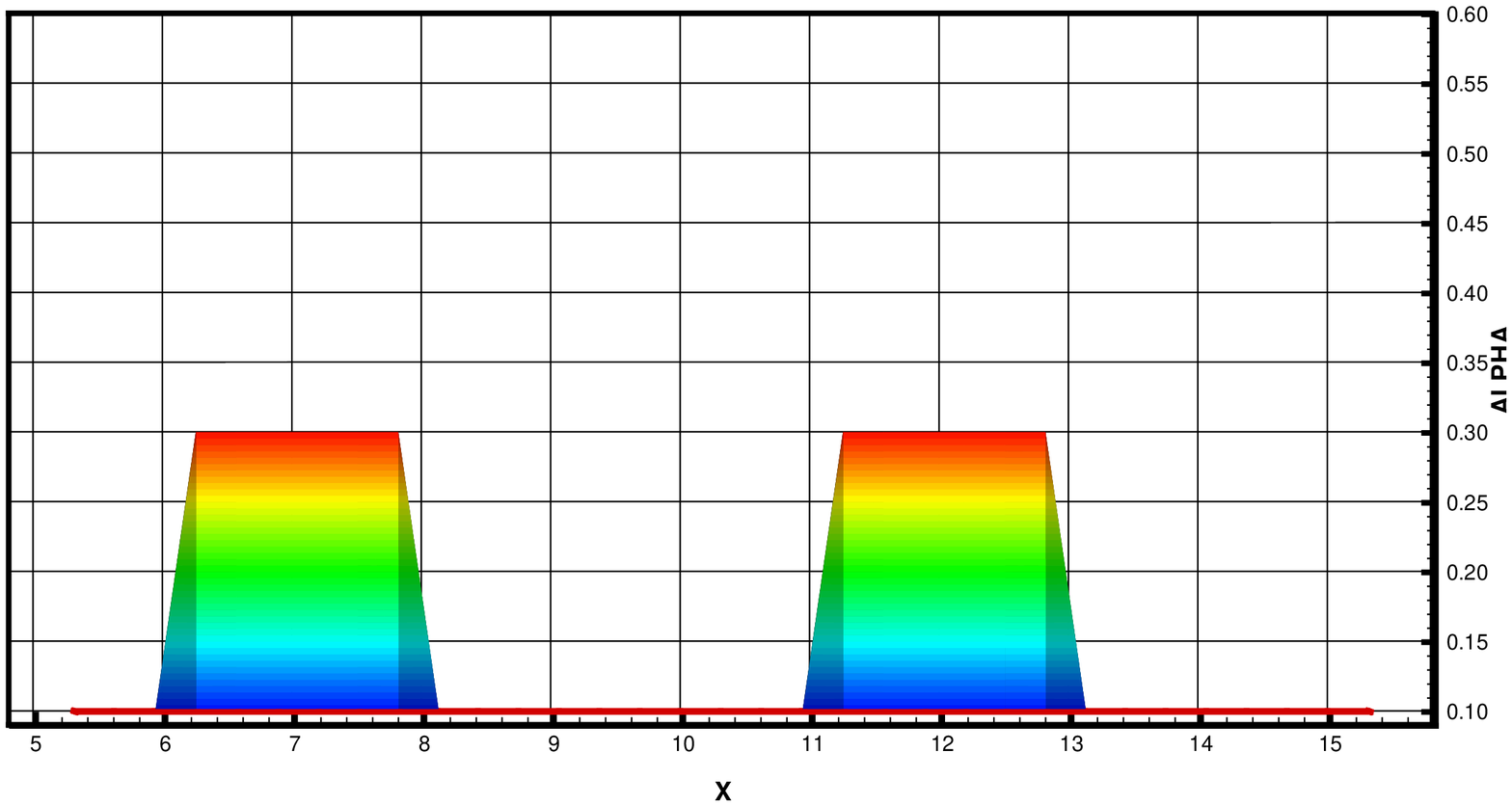} 
\end{figure}
\textit{Figure 5a}

\begin{figure}[!h]
\begin{picture}(0,0)\end{picture}
\includegraphics[width=2in,height=2in,angle=0]{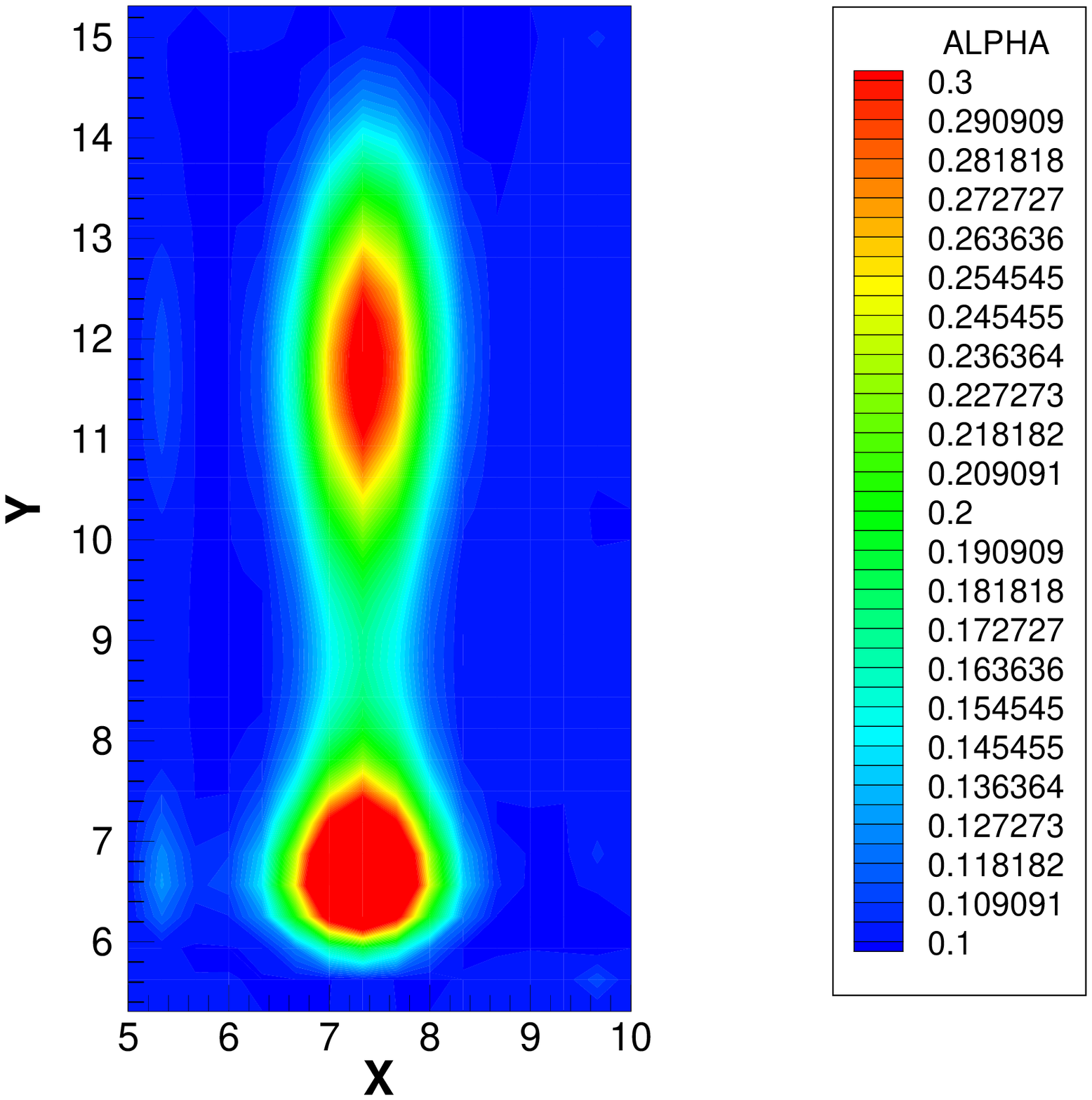} \setlength{%
\unitlength}{1in} 
\hspace{0.5in} \begin{picture}(0,0)\end{picture}
\includegraphics[width=2in,height=2in,angle=0]{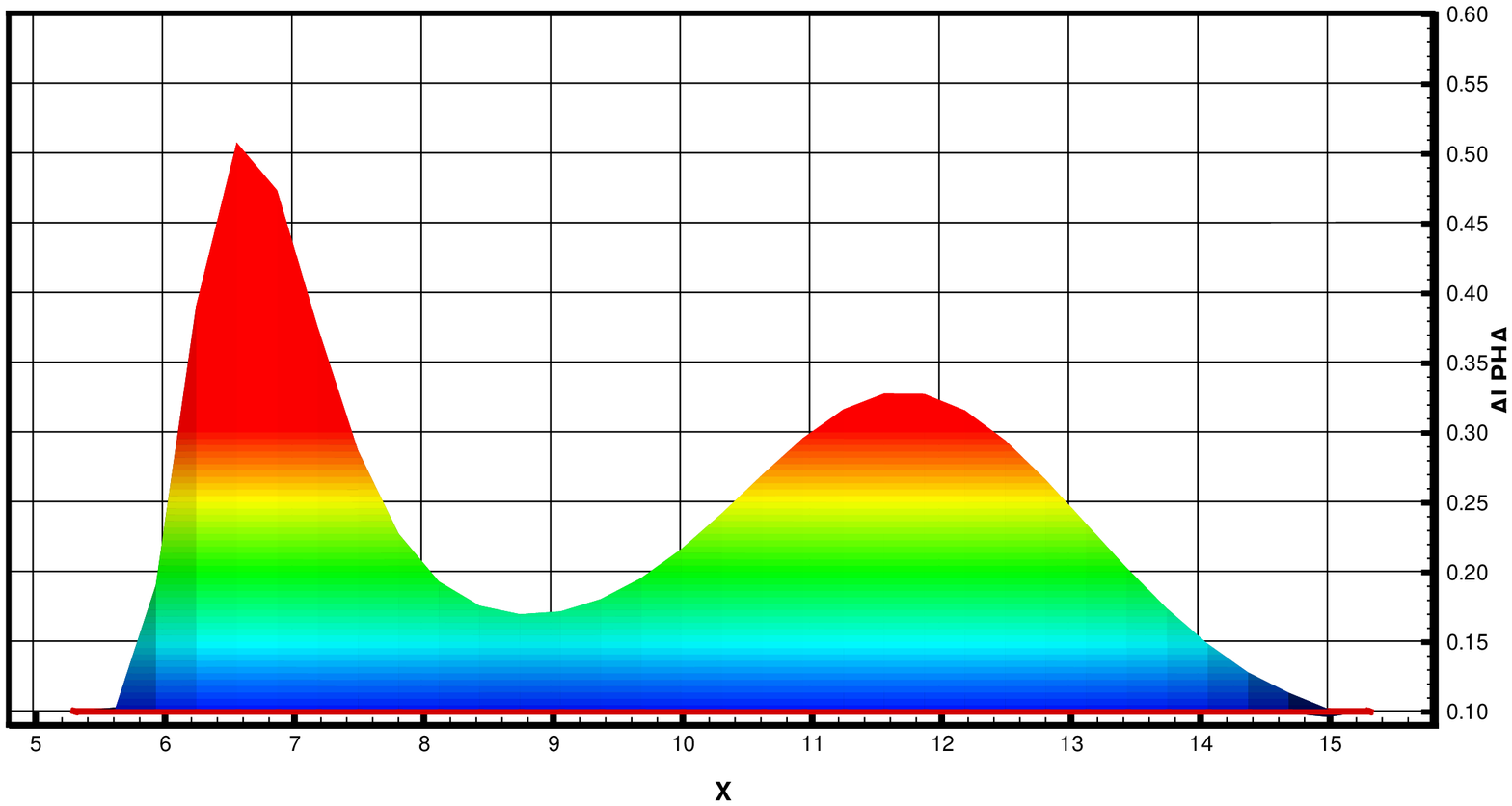} 
\end{figure}
\textit{Figure 5b}

\begin{figure}[!h]
\begin{picture}(0,0)\end{picture}
\includegraphics[width=2in,height=2in,angle=0]{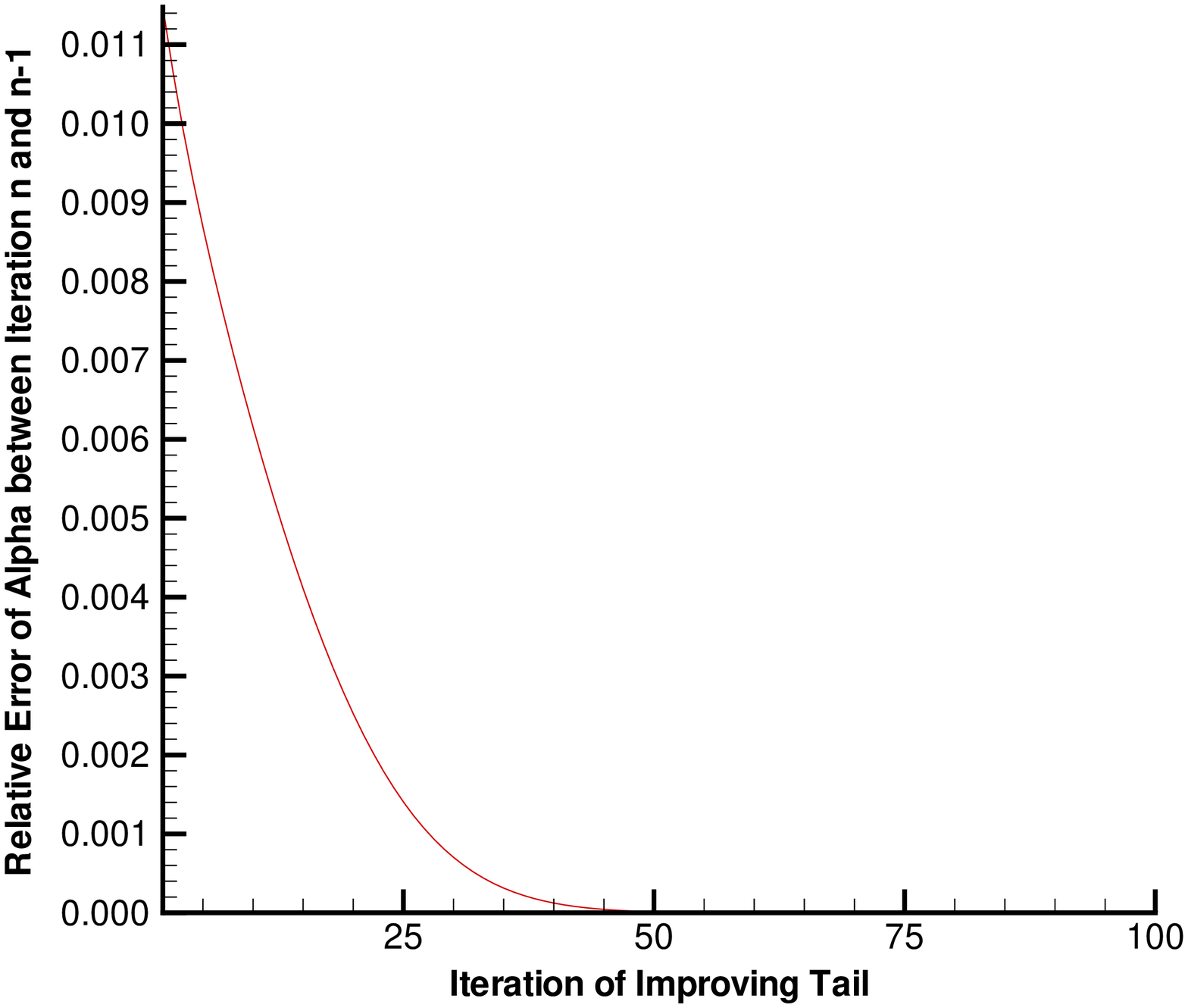} \setlength{%
\unitlength}{1in} 
\hspace{0.5in} \begin{picture}(0,0)\end{picture}
\includegraphics[width=2in,height=2in,angle=0]{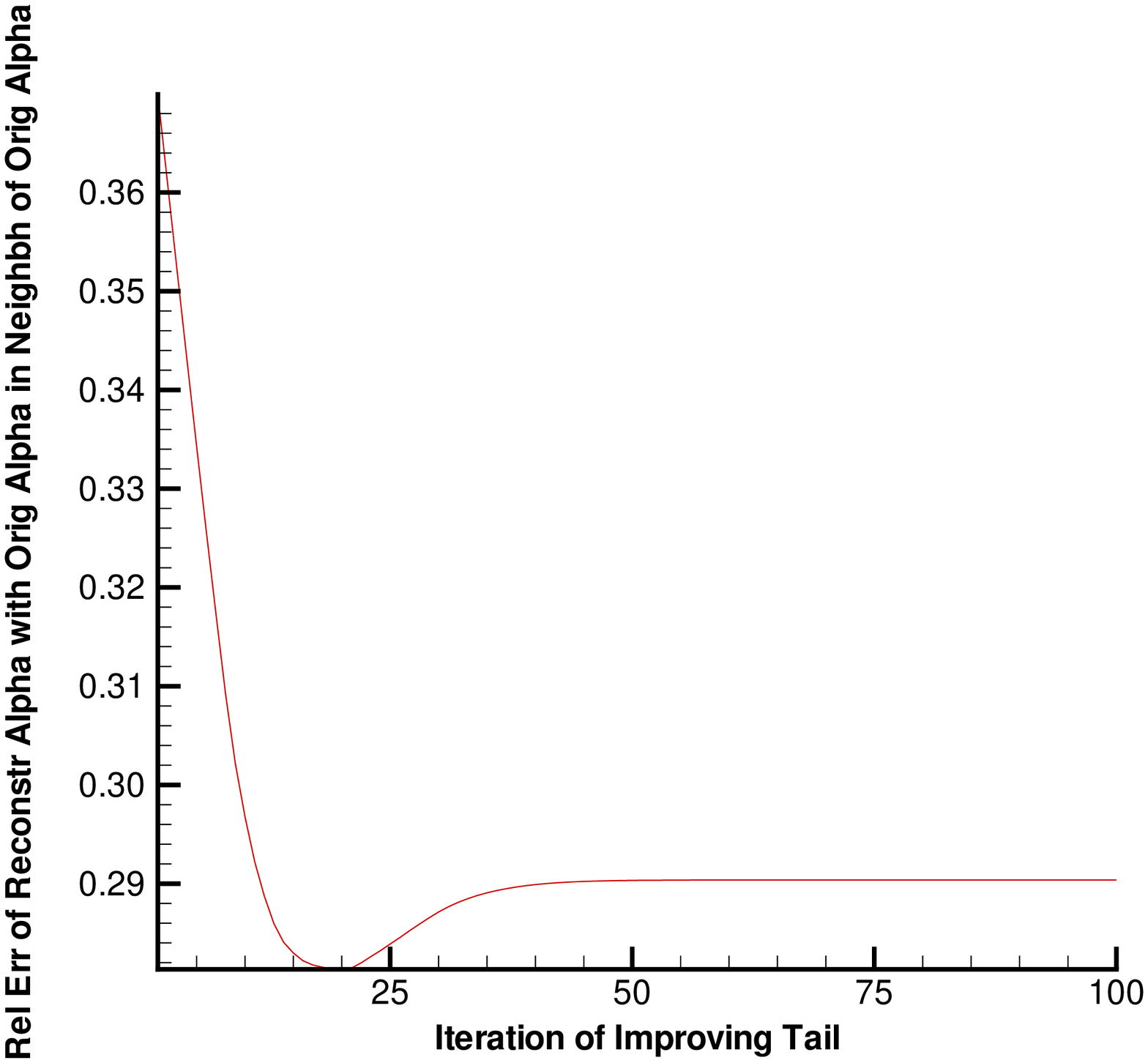} 
\end{figure}
\textit{Figure 5c on left and Figure 5d on right}

\textit{Figures 5a-5d. Figure 5a displays the original distribution. Figure
5b shows reconstruction from the noisy data, the level of noise is $2\%$. The 1-d cross sections next to
the 2-d figures show the profiles of the original inclusion and its
reconstruction at $z=7$ }

\textit{We illustrate on Figure 5c the difference of two consecutive
reconstruction 
\begin{equation*}
||a_{m}(x,z)-a_{m-1}(x,z)||\equiv {\frac{\sqrt{\sum%
\limits_{i=1,...,i_{max},j=1,...,j_{max}}|(a_{m}(x_{i},z_{j})-a_{m-1}(x_{i},z_{j}))|^{2}%
}}{\sqrt{N_{1}}\max {|a_{m-1}(x_{i},z_{j})|}}}
\end{equation*}
as a function of the number of iteration $m$. The function is used for
determining stopping criterion. Figure 5d depicts the relative error in
comparison with actual inclusion expressed by 
\begin{equation*}
RMSE\equiv {\frac{\sqrt{\sum%
\limits_{i=1,...,i_{max},j=1,...,j_{max}}|(a_{m}(x_{i},z_{j})-a(x_{i},z_{j}))|^{2}%
}}{\sqrt{N_{1}}\max {|a(x_{i},z_{j})|}}}
\end{equation*}
as a function of the number of iterations $m$. }

\newpage

\begin{figure}[!h]
\begin{picture}(0,0)\end{picture}
\includegraphics[width=2in,height=2in,angle=0]{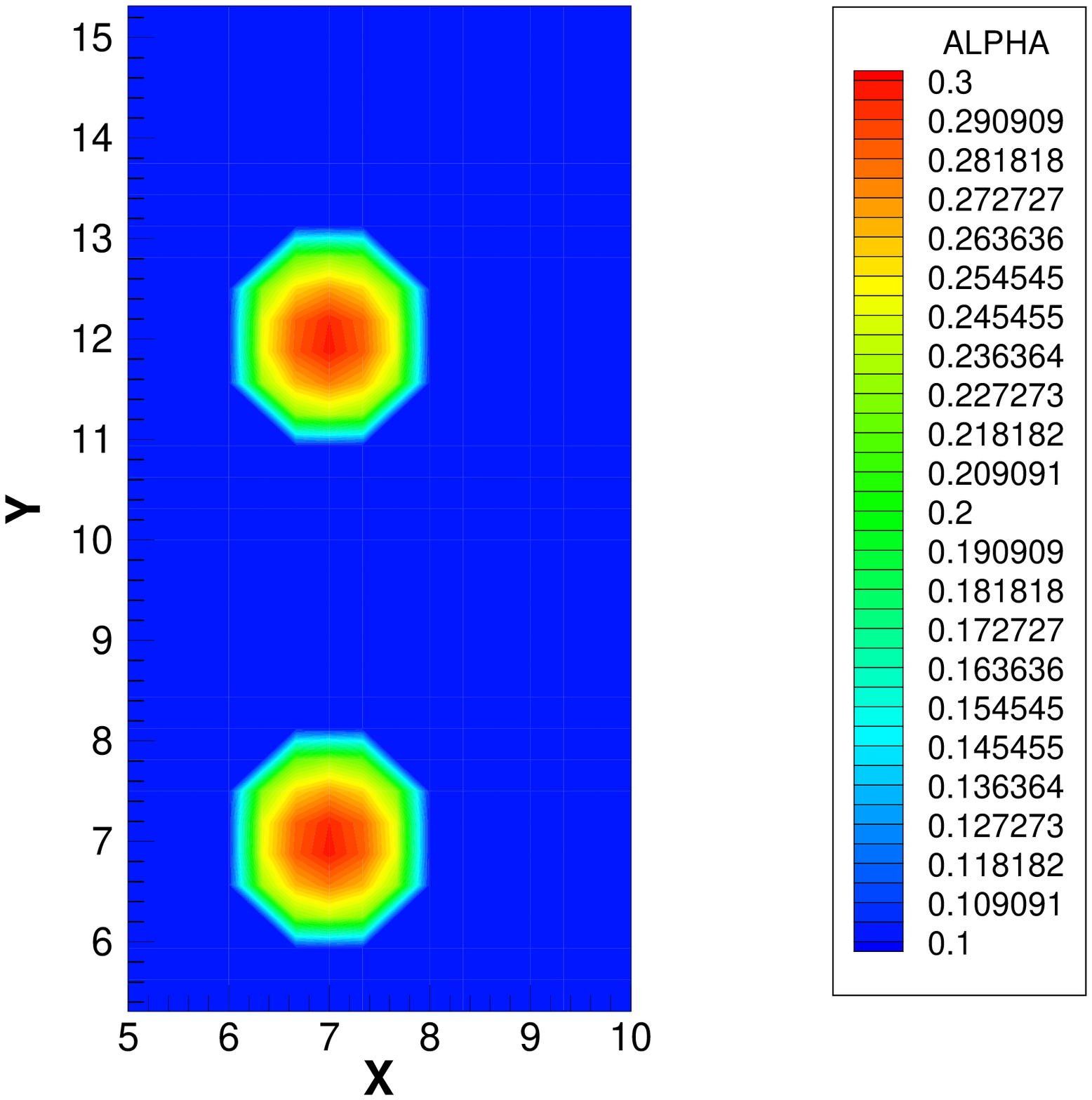} \setlength{%
\unitlength}{1in} 
\hspace{0.5in} \begin{picture}(0,0)\end{picture}
\includegraphics[width=2in,height=2in,angle=0]{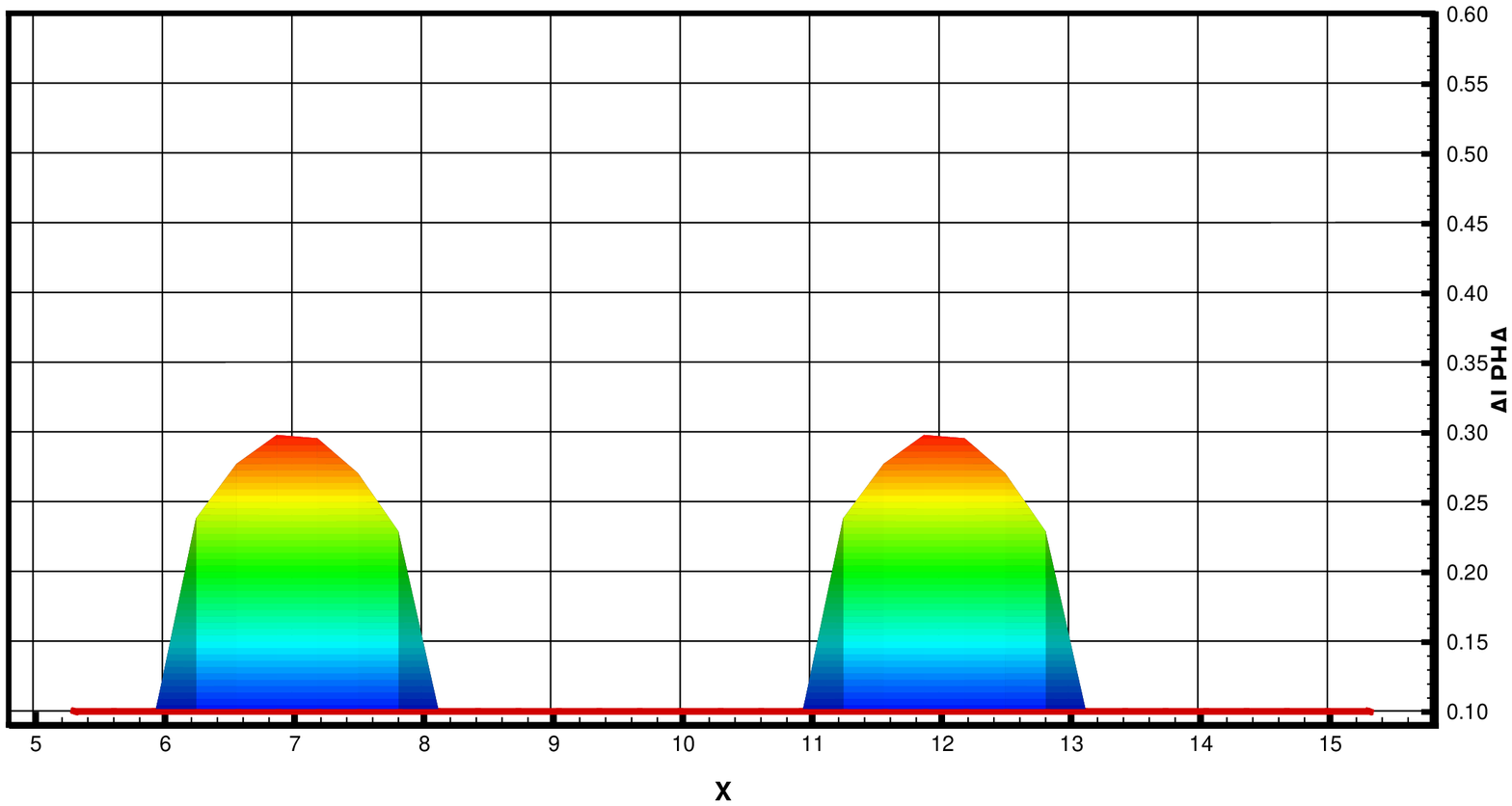} 
\end{figure}
\textit{Figure 6a}

\begin{figure}[!h]
\begin{picture}(0,0)\end{picture}
\includegraphics[width=2in,height=2in,angle=0]{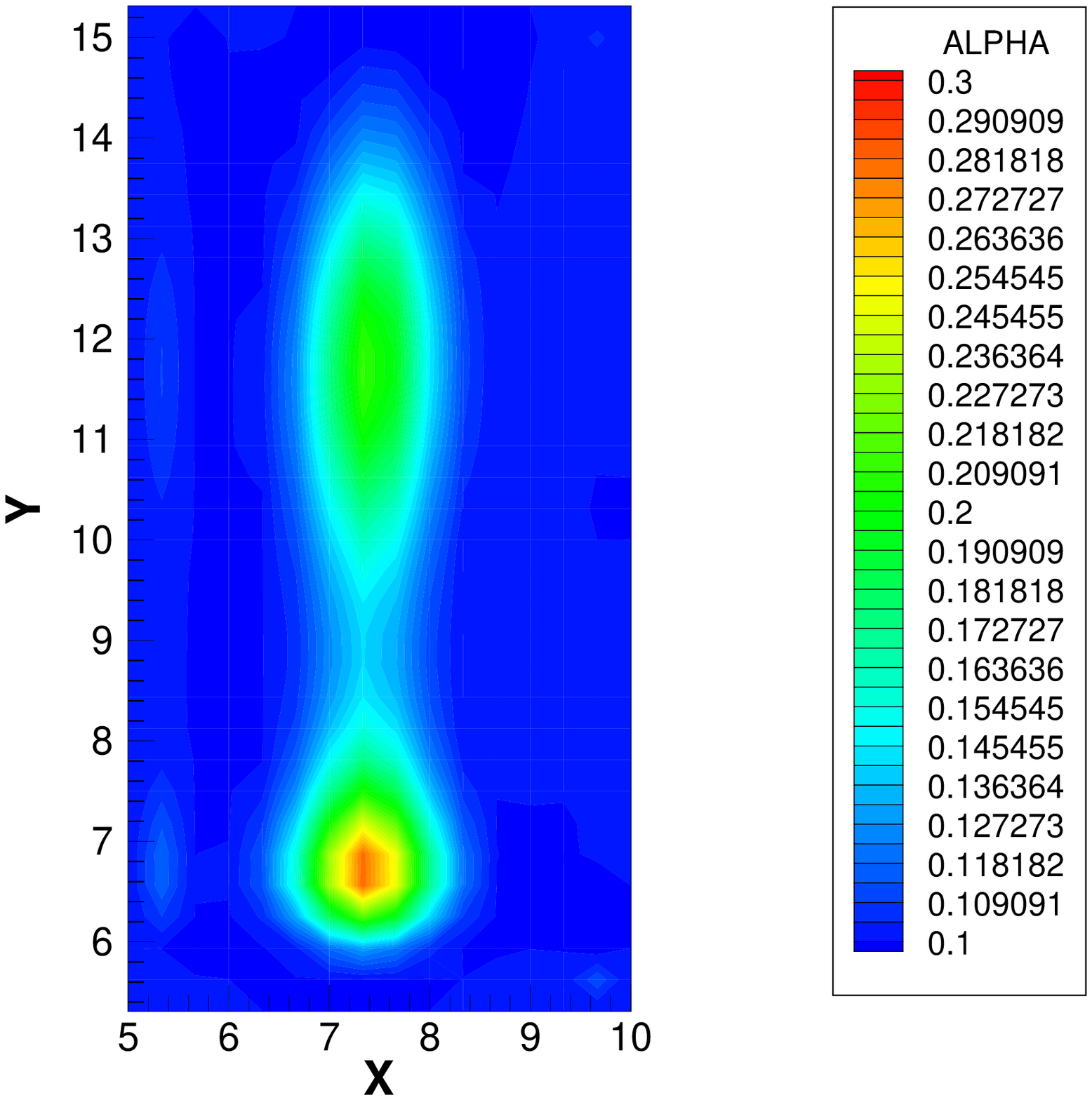} \setlength{%
\unitlength}{1in} 
\hspace{0.5in} \begin{picture}(0,0)\end{picture}
\includegraphics[width=2in,height=2in,angle=0]{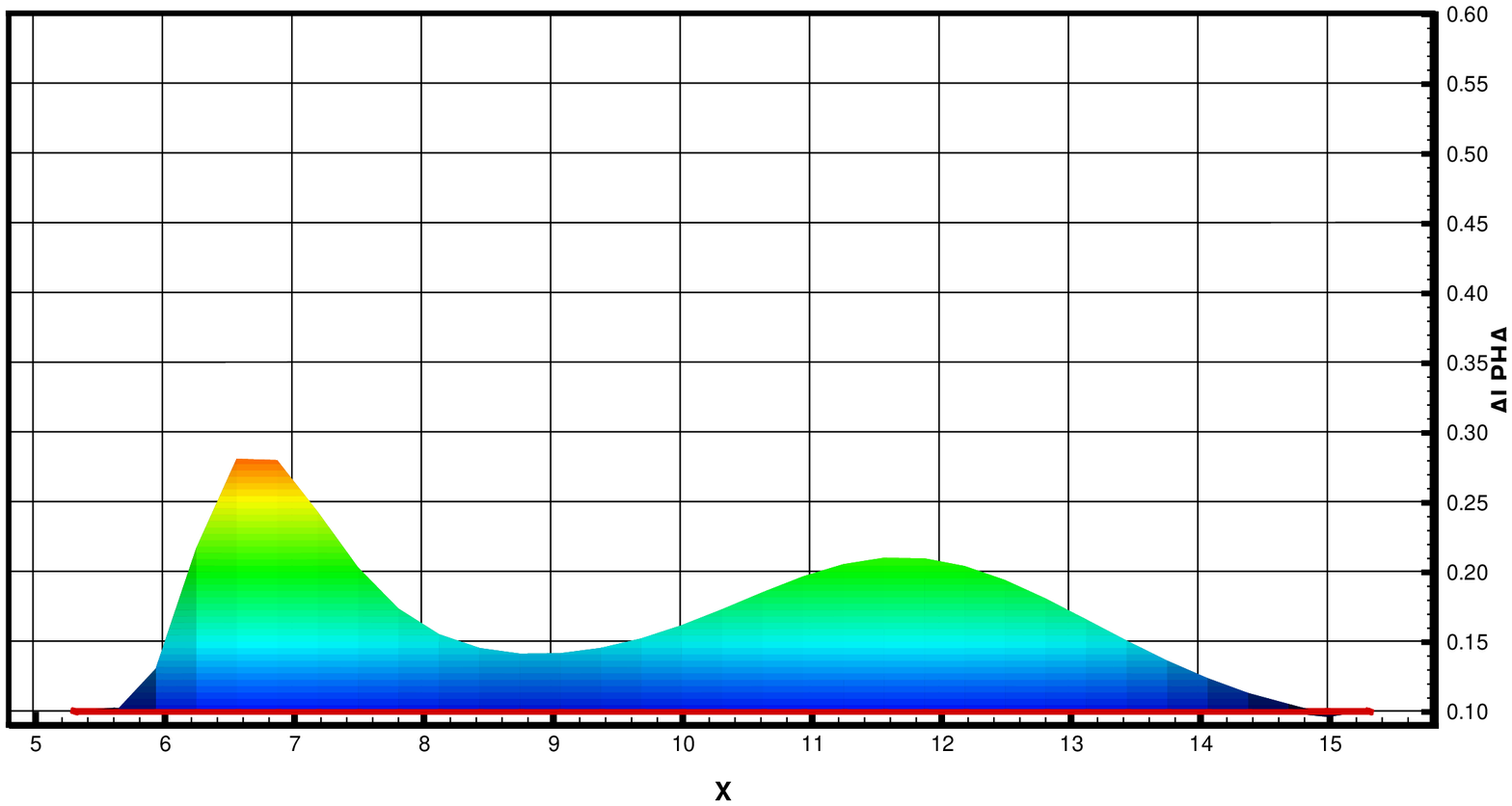} 
\end{figure}
\textit{Figure 6b}

\begin{figure}[!h]
\begin{picture}(0,0)\end{picture}
\includegraphics[width=2in,height=2in,angle=0]{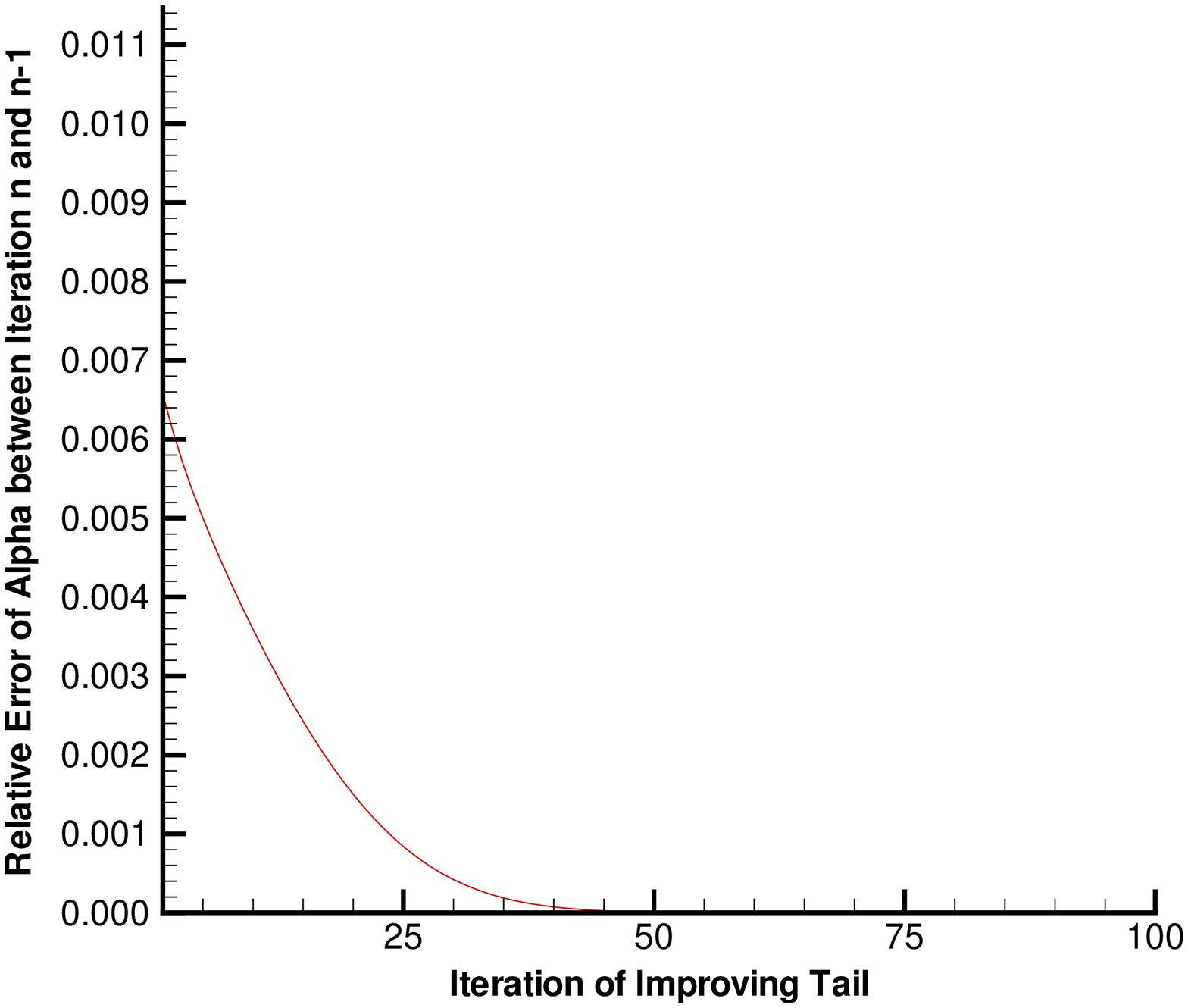} \setlength{%
\unitlength}{1in} 
\hspace{0.5in} \begin{picture}(0,0)\end{picture}
\includegraphics[width=2in,height=2in,angle=0]{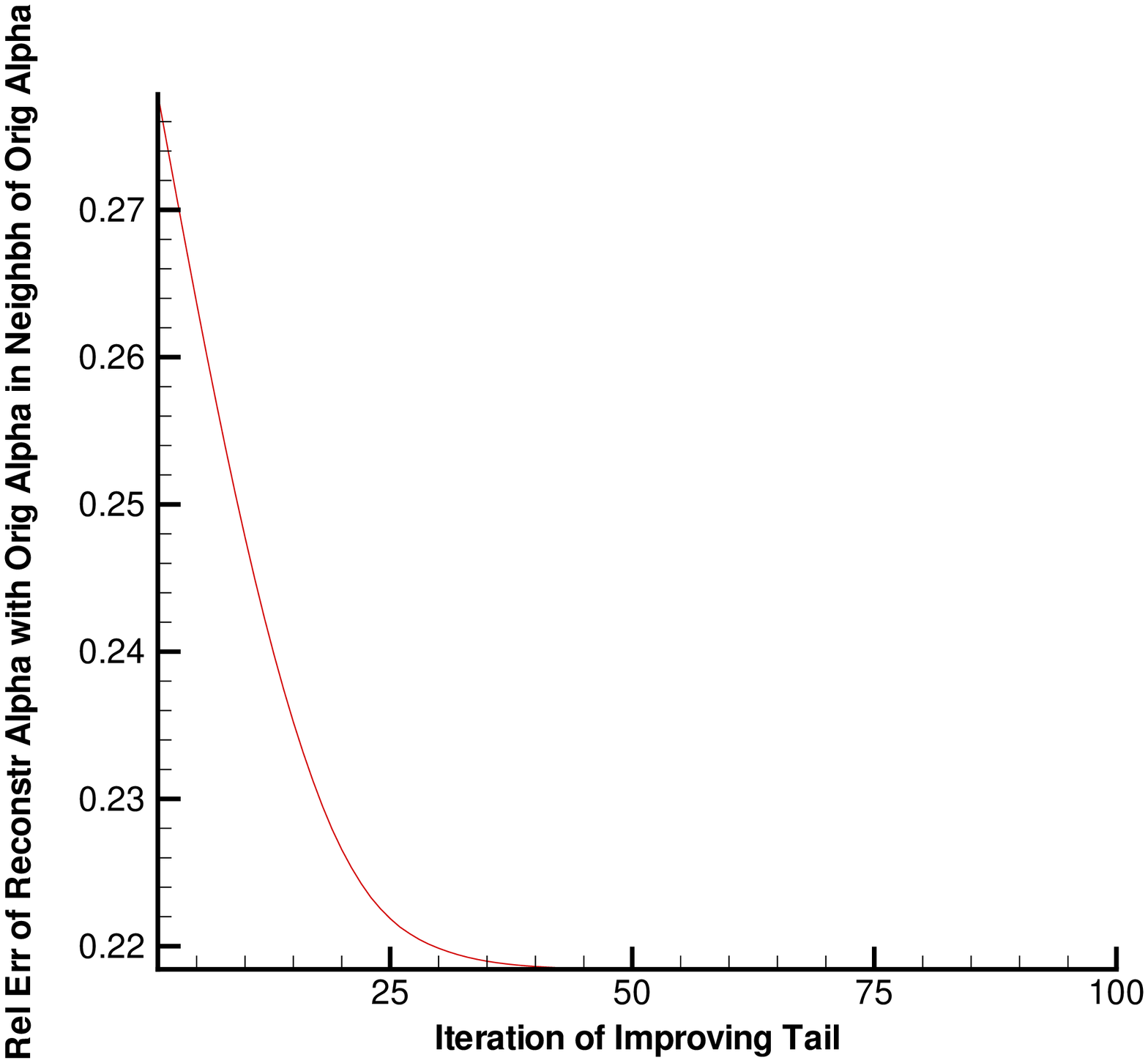} 
\end{figure}
\textit{Figure 6c on left and Figure 6d on right}

\textit{Figures 6a-6d. Figure 6a displays the original function in \ two
inclusions. Figure 6b shows the reconstruction result with 2\% noise.} 
\textit{We illustrate on Figure 6c the difference of two consecutive
reconstruction 
\begin{equation*}
||a_{m}(x,z)-a_{m-1}(x,z)||\equiv {\frac{\sqrt{\sum%
\limits_{i=1,...,i_{max},j=1,...,j_{max}}|(a_{m}(x_{i},z_{j})-a_{m-1}(x_{i},z_{j}))|^{2}%
}}{\sqrt{N_{1}}\max {|a_{m-1}(x_{i},z_{j})|}}}
\end{equation*}
as a function of the number of iteration $m$. The function is used for
determining stopping criterion. Figure 6d depicts the relative error in
comparison with actual inclusion expressed by 
\begin{equation*}
RMSE\equiv {\frac{\sqrt{\sum%
\limits_{i=1,...,i_{max},j=1,...,j_{max}}|(a_{m}(x_{i},z_{j})-a(x_{i},z_{j}))|^{2}%
}}{\sqrt{N_{1}}\max {|a(x_{i},z_{j})|}}}
\end{equation*}
as a function of the number of iterations $m$.}

\newpage \vskip5in

\begin{figure}[!h]
\begin{picture}(0,0)\end{picture}
\includegraphics[width=2in,height=2in,angle=0]{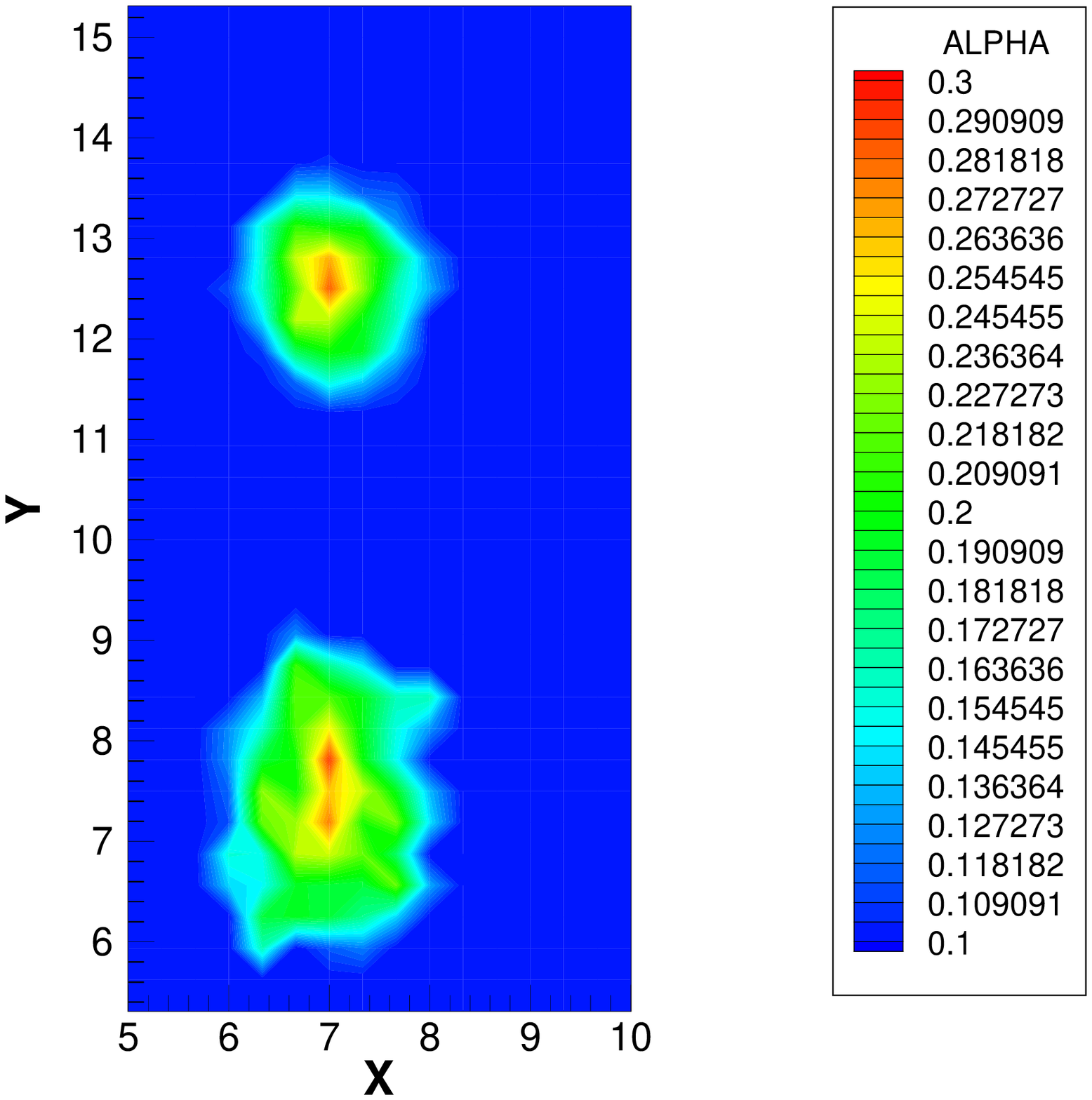} \setlength{%
\unitlength}{1in} 
\hspace{0.5in} \begin{picture}(0,0)\end{picture}
\includegraphics[width=2in,height=2in,angle=0]{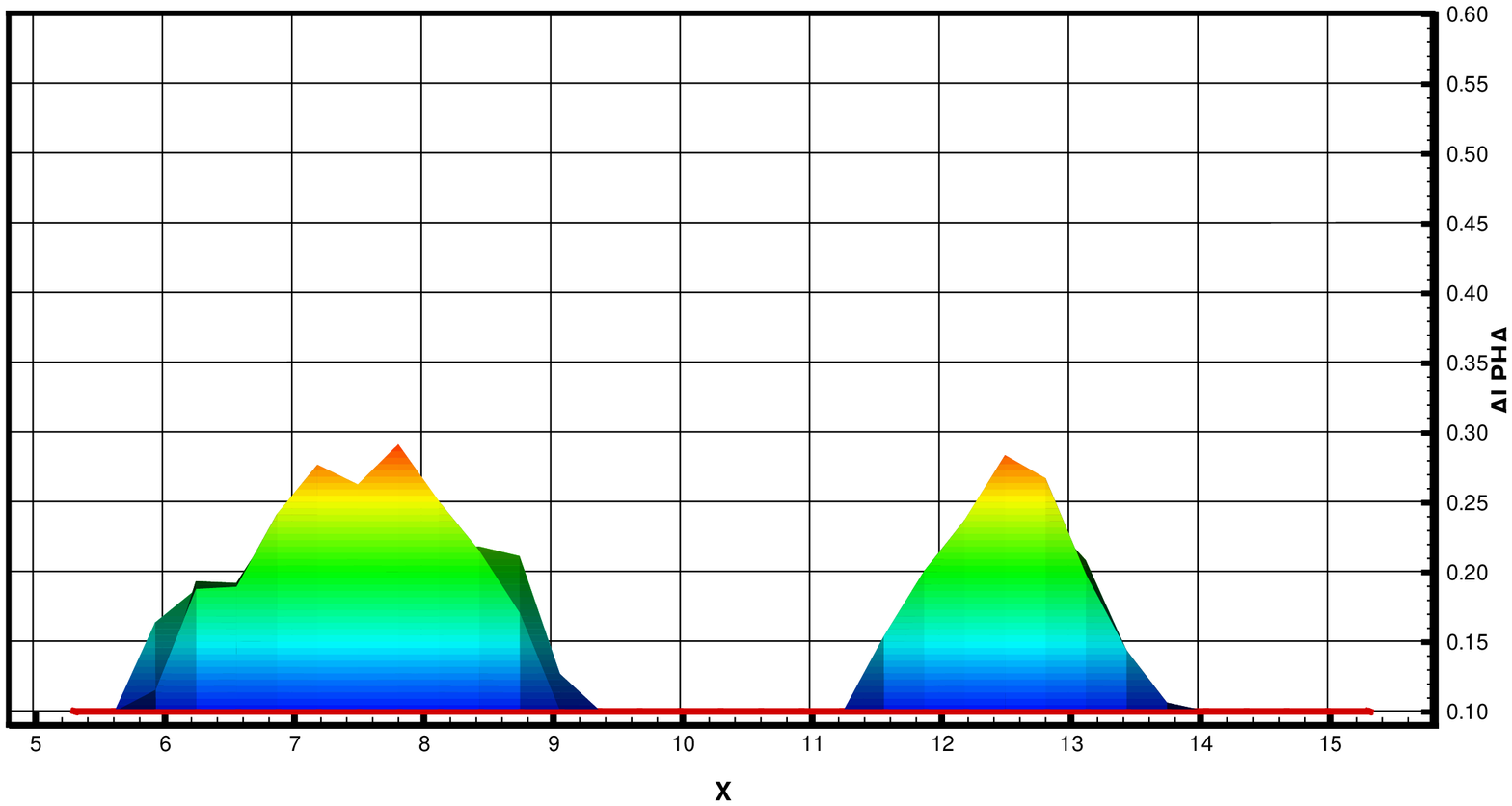} 
\end{figure}
\textit{Figure 7a}

\begin{figure}[!h]
\begin{picture}(0,0)\end{picture}
\includegraphics[width=2in,height=2in,angle=0]{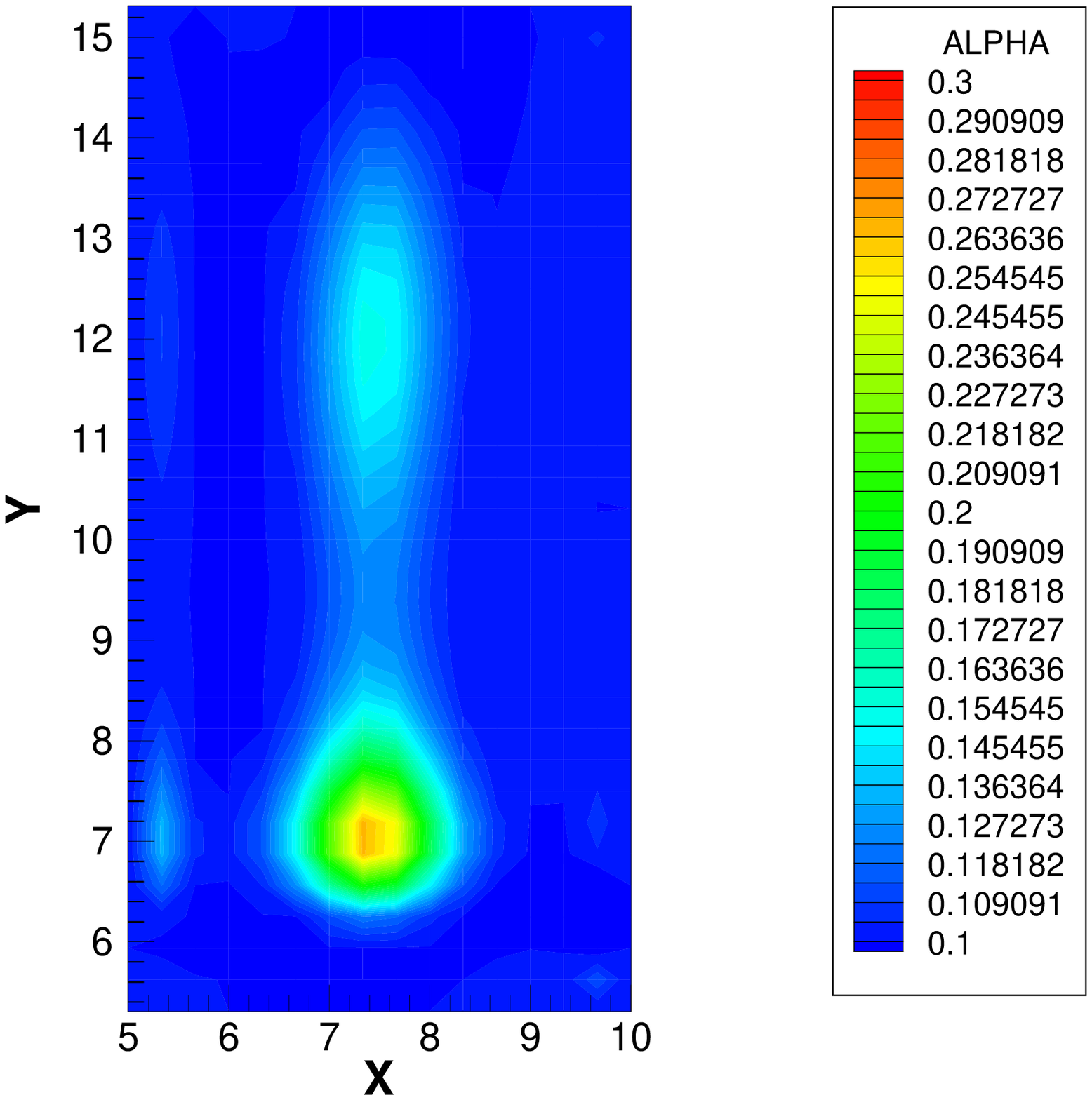} \setlength{%
\unitlength}{1in} 
\hspace{0.5in} \begin{picture}(0,0)\end{picture}
\includegraphics[width=2in,height=2in,angle=0]{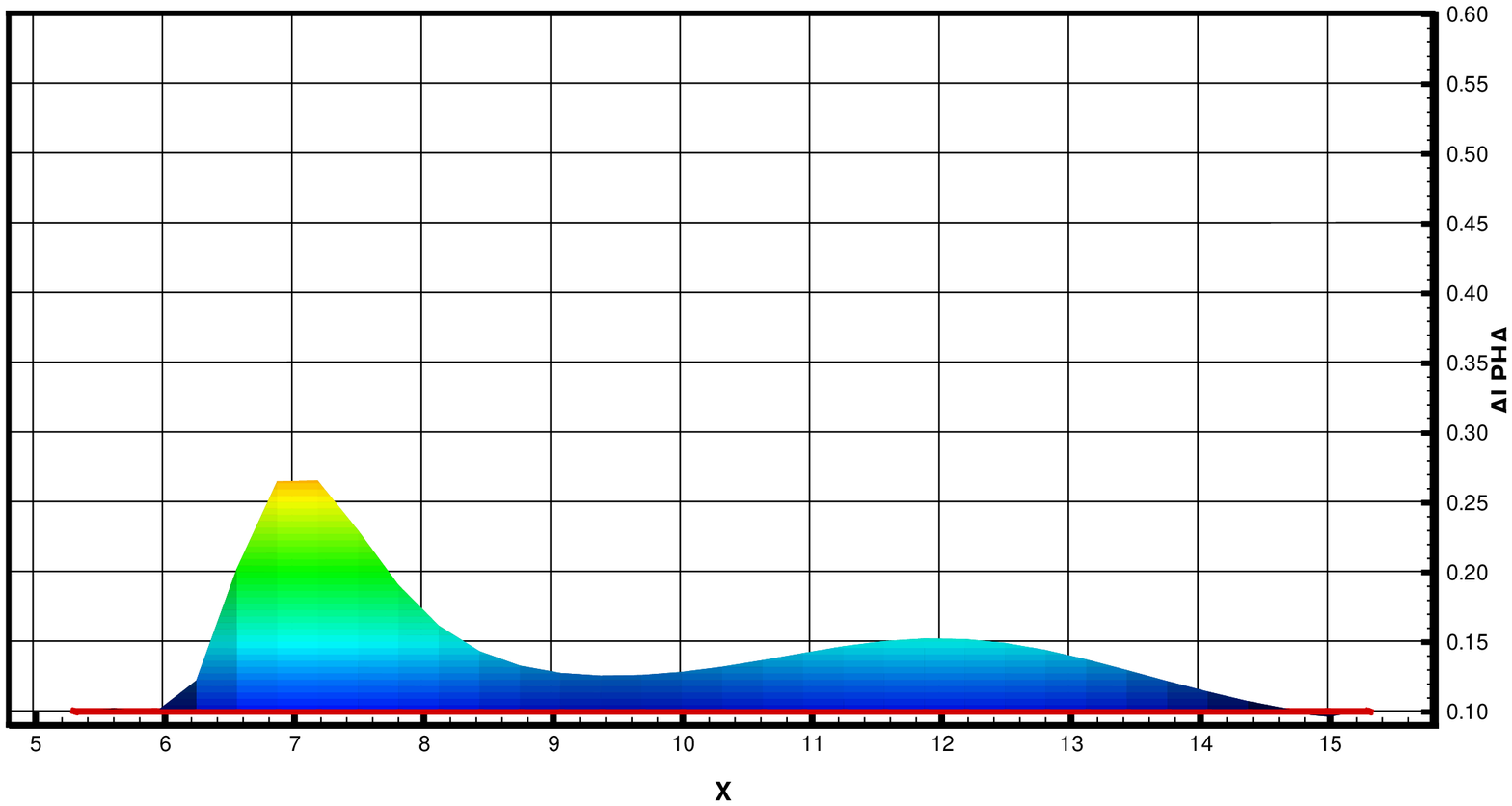} 
\end{figure}
\textit{Figure 7b}

\begin{figure}[!h]
\begin{picture}(0,0)\end{picture}
\includegraphics[width=2in,height=2in,angle=0]{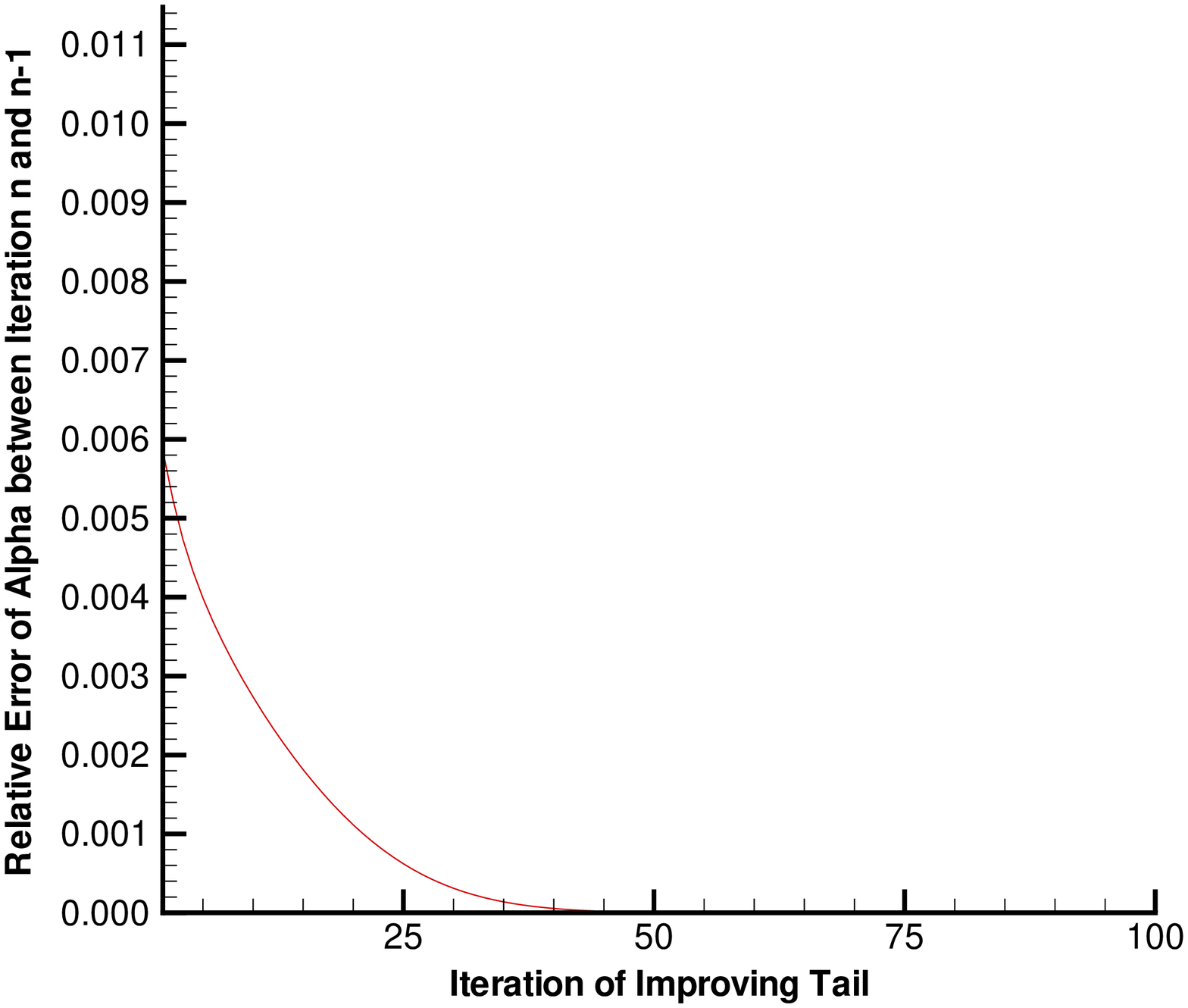} \setlength{%
\unitlength}{1in} 
\hspace{0.5in} \begin{picture}(0,0)\end{picture}
\includegraphics[width=2in,height=2in,angle=0]{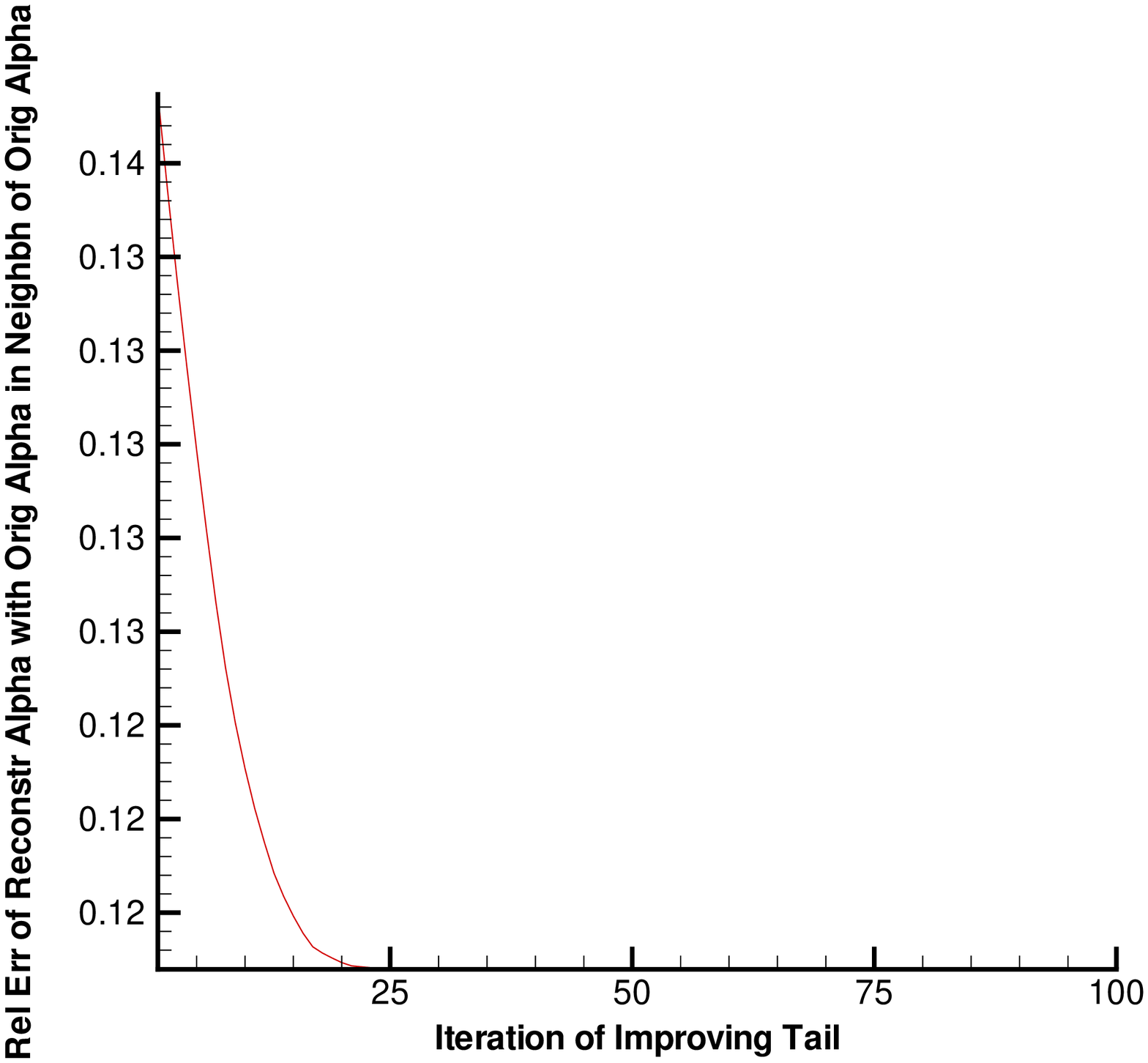} 
\end{figure}
\textit{Figure 7c on left and Figure 7d on right}

\textit{Figures 7a-7d. Figure 6a displays the original distribution. Figure
7b displays the reconstruction result with 2\% noise in the data. The 1-d
cross sections show the profiles of the original inclusion and its
reconstruction at $z=7$ } \textit{We illustrate on Figure 7c the difference
of two consecutive reconstruction 
\begin{equation*}
||a_{m}(x,z)-a_{m-1}(x,z)||\equiv {\frac{\sqrt{\sum%
\limits_{i=1,...,i_{max},j=1,...,j_{max}}|(a_{m}(x_{i},z_{j})-a_{m-1}(x_{i},z_{j}))|^{2}%
}}{\sqrt{N_{1}}\max {|a_{m-1}(x_{i},z_{j})|}}}
\end{equation*}
as a function of the number of iteration $m$. The function is used for
determining stopping criterion. Figure 7d depicts the relative error in
comparison with actual inclusion expressed by 
\begin{equation*}
RMSE\equiv {\frac{\sqrt{\sum%
\limits_{i=1,...,i_{max},j=1,...,j_{max}}|(a_{m}(x_{i},z_{j})-a(x_{i},z_{j}))|^{2}%
}}{\sqrt{N_{1}}\max {|a(x_{i},z_{j})|}}}
\end{equation*}
as a function of the number of iterations $m$.}

\end{document}